\documentclass[reqno, 11pt, a4paper]{amsart}

\makeatletter
\renewcommand{\author}[2][]{%
  \ifx\@empty\authors
    \gdef\authors{#2}%
  \else
    \g@addto@macro\authors{\and#2}%
    \g@addto@macro\addresses{\author{}}%
  \fi
  \@ifnotempty{#1}{%
    \ifx\@empty\shortauthors
      \gdef\shortauthors{\footnotesize\scshape#1}%
    \else
      \g@addto@macro\shortauthors{\and\footnotesize\scshape#1}%
    \fi
  }}
\edef\author{\@nx\@dblarg
  \@xp\@nx\csname\string\author\endcsname}
\let\shortauthors\@empty   \let\authors\@empty
\def\maketitle{\par
  \@topnum\z@%
  \@setcopyright
  \thispagestyle{firstpage}%
  \uppercasenonmath\shorttitle
  \ifx\@empty\shortauthors \let\shortauthors\shorttitle
  \else \andify\shortauthors
  \fi
  \@maketitle@hook
  \begingroup
  \@maketitle
  \toks@\@xp{\shortauthors}\@temptokena\@xp{\shorttitle}%
  \toks4{\def\\{ \ignorespaces}}%
  \edef\@tempa{%
    \@nx\markboth{\the\toks4
      \@nx{\the\toks@}}{\the\@temptokena}}%
  \@tempa
  \endgroup
  \c@footnote\z@
  \@cleartopmattertags}
\def\@setauthors{%
  \begingroup
  \def\thanks{\protect\thanks@warning}%
  \trivlist
  \centering\@topsep30\p@\relax
  \advance\@topsep by -\baselineskip
  \item\relax
  \author@andify\authors
  \def\\{\protect\linebreak}%
  \small\scshape\authors%
  \ifx\@empty\contribs
  \else
    ,\penalty-3 \space \@setcontribs
    \@closetoccontribs
  \fi
  \endtrivlist
  \endgroup}
\renewenvironment{abstract}{%
  \ifx\maketitle\relax
    \ClassWarning{\@classname}{Abstract should precede
      \protect\maketitle\space in AMS document classes; reported}%
  \fi
  \global\setbox\abstractbox=\vtop \bgroup
    \normalfont\Small
    \list{}{\labelwidth\z@
      \leftmargin3pc \rightmargin\leftmargin
      \listparindent\normalparindent \itemindent\z@
      \parsep\z@ \@plus\p@
      }%
    \item[\hskip\labelsep\bfseries\abstractname.]%
}{%
  \endlist\egroup
  \ifx\@setabstract\relax \@setabstracta \fi}

\renewcommand{\tocsection}[3]{%
  \indentlabel{\@ifnotempty{#2}{\ignorespaces#1 #2. \ }}#3}
\renewcommand\datename{Date}
\def\@setdate{{\itshape \datename.}\enspace\@date\@addpunct.}
\renewcommand{\@secnumfont}{\bfseries}
\def\section{\@startsection{section}{1}%
  \z@{.7\linespacing\@plus\linespacing}{.5\linespacing}%
  {\normalfont\bfseries\centering}}
\def\subsection{\@startsection{subsection}{2}%
  \z@{.5\linespacing\@plus.7\linespacing}{-.5em}%
  {\normalfont\itshape}}
\def\@captionheadfont{\bfseries}
\def\curraddrname{{\itshape Current address.}}
\def\emailaddrname{{\itshape Email address.}}
\def\urladdrname{{\itshape URL.}}
\def\@setaddresses{\par
  \nobreak \begingroup
  \footnotesize
  \def\author##1{\nobreak\addvspace\bigskipamount}%
  \def\\{\unskip, \ignorespaces}%
  \interlinepenalty\@M
  \def\address##1##2{\begingroup
    \par\addvspace\bigskipamount\indent
    \@ifnotempty{##1}{(\ignorespaces##1\unskip) }%
    {\ignorespaces##2}\par\endgroup}%
  \def\curraddr##1##2{\begingroup
    \@ifnotempty{##2}{\nobreak\indent\curraddrname
      \@ifnotempty{##1}{, \ignorespaces##1\unskip}\enspace
      ##2\par}\endgroup}%
  \def\email##1##2{\begingroup
    \@ifnotempty{##2}{\nobreak\indent\emailaddrname
      \@ifnotempty{##1}{, \ignorespaces##1\unskip}\enspace
      \ttfamily##2\par}\endgroup}%
  \def\urladdr##1##2{\begingroup
    \def~{\char`\~}%
    \@ifnotempty{##2}{\nobreak\indent\urladdrname
      \@ifnotempty{##1}{, \ignorespaces##1\unskip}\enspace
      \ttfamily##2\par}\endgroup}%
  \addresses
  \endgroup}
\makeatother

\usepackage[utf8]{inputenc}
\usepackage[T1]{fontenc}
\usepackage{lmodern}
\usepackage{amssymb}
\usepackage[colorlinks=true,citecolor=blue,urlcolor=blue,linkcolor=blue,backref=page,ocgcolorlinks]{hyperref}
\usepackage{graphicx}
\usepackage{tikz-cd}
\usepackage{mathrsfs}
\usepackage{mathtools}

\numberwithin{equation}{section}
\def\secskip{\vspace{.5\linespacing plus.7\linespacing}}
\newcommand\subsectionempty{\subsection{\texorpdfstring{\!\!\!}{}}}

\renewcommand\backrefalt[4]{%
  \ifcase #1%
    (Not cited)%
  \or%
    (Page~#2)%
  \else%
    (Pages~#2)%
  \fi}

\newcommand*{\bP}{\mathbb{P}}
\newcommand*{\bQ}{\mathbb{Q}}

\newcommand*{\bZ}{\mathbb{Z}}

\newcommand*{\sA}{\mathscr{A}}

\newcommand*{\sC}{\mathscr{C}}
\newcommand*{\sD}{\mathscr{D}}

\newcommand*{\sF}{\mathscr{F}}
\newcommand*{\sG}{\mathscr{G}}

\newcommand*{\sJ}{\mathscr{J}}

\newcommand*{\sL}{\mathscr{L}}
\newcommand*{\sM}{\mathscr{M}}

\newcommand*{\sO}{\mathscr{O}}
\newcommand*{\sP}{\mathscr{P}}

\newcommand*{\sR}{\mathscr{R}}

\newcommand*{\sT}{\mathscr{T}}

\newcommand*{\rd}{\mathrm{d}}

\newcommand*{\fS}{\mathfrak{S}}

\renewcommand*{\epsilon}{\varepsilon}

\newcommand*{\apriori}{\textit{a priori} }

\newcommand*{\etc}{\textit{etc.}}
\newcommand*{\ie}{\textit{i.e.}~}
\newcommand*{\loccit}{\textit{loc.\,cit.}}
\newcommand*{\resp}{resp.~}

\newcommand*{\rand}{\textup{ \ and \ }}
\newcommand*{\rfor}{\textup{ \ for \ }}
\newcommand*{\rforall}{\textup{ \ for all \ }}

\newcommand*{\rif}{\textup{ \ if \ }}

\newcommand*{\rin}{\textup{ \ in \ }}
\newcommand*{\ron}{\textup{ \ on \ }}

\newcommand*{\rwith}{\textup{ \ with \ }}

\newcommand*{\rleft}{\textup{(}}
\newcommand*{\rright}{\textup{)}}

\newcommand*{\ct}{\mathop\mathrm{ct}\nolimits}
\newcommand*{\CH}{\mathop\mathrm{CH}\nolimits}

\newcommand*{\End}{\mathop\mathrm{End}\nolimits}
\newcommand*{\ev}{\mathop\mathrm{ev}\nolimits}

\newcommand*{\id}{\mathop\mathrm{id}\nolimits}

\newcommand*{\mon}{\mathop\mathrm{mon}\nolimits}
\newcommand*{\Mon}{\mathop\mathrm{Mon}\nolimits}

\newcommand*{\Pic}{\mathop\mathrm{Pic}\nolimits}
\newcommand*{\pr}{\mathop\mathrm{pr}\nolimits}

\newcommand*{\prim}{\mathop\mathrm{prim}\nolimits}
\newcommand*{\rt}{\mathop\mathrm{rt}\nolimits}

\newcommand*{\SL}{\mathop\mathrm{SL}\nolimits}

\newcommand*{\fsl}{\mathop\mathfrak{sl}\nolimits}

\DeclareRobustCommand{\longtwoheadrightarrow}{\relbar\!\!\!\joinrel\twoheadrightarrow}

\newcommand*{\into}{\hookrightarrow}
\newcommand*{\onto}{\longtwoheadrightarrow}

\newcommand*{\xto}{\xrightarrow}

\makeatletter
\newcommand*{\vast}{\bBigg@{3}}
\newcommand*{\Vast}{\bBigg@{4}}
\makeatother

\begin{document}

\title[A tale of two tautological rings]{Cycles on curves and Jacobians:\\a tale of two tautological rings}
\author{Qizheng Yin}
\date{\today}
\address{ETH Zürich, Departement Mathematik, Rämistrasse 101, 8092 Zürich, Switzerland}
\email{\href{mailto:yin@math.ethz.ch}{yin@math.ethz.ch}}
\thanks{The author is supported by grant ERC-2012-AdG-320368-MCSK}

\begin{abstract}
We connect two notions of tautological ring: one for the moduli space of curves (after Mumford, Faber, \etc), and the other for the Jacobian of a curve (after Beauville, \mbox{Polishchuk}, \etc). The motivic Lefschetz decomposition on the Jacobian side produces relations between tautological classes, leading to results about Faber's Gorenstein conjecture on the curve side. We also relate certain Gorenstein properties on both sides and verify them for small genera. Further, we raise the question whether all tautological relations are motivic, giving a possible explanation why the Gorenstein properties may not hold.
\end{abstract}
\maketitle

\tableofcontents

\section*{Introduction}
\secskip

\noindent Tautological classes are classes of algebraic cycles carrying certain geometric information. They form a ring called the tautological ring. Classically there are two such notions of tautological ring: one introduced by Mumford \cite{Mum83} for the moduli space of curves, and the other by Beauville \cite{Bea04} for the Jacobian of a curve (and modulo algebraic equivalence). 

\subsection{Curve side} --- Denote by $\sM_g$ the moduli space of smooth curves of genus $g$, and by $p \colon \sC_g \to \sM_g$ the universal curve. Consider the classes $\kappa_i = p_*(K^{i + 1})$, with $K$ the relative canonical divisor. The tautological ring $\sR(\sM_g)$ is defined to be the $\bQ$-subalgebra of the Chow ring $\CH(\sM_g)$ (with $\bQ$-coeffcients) generated by $\{\kappa_i\}$.

The study of $\sR(\sM_g)$ is partly motivated by Faber's Gorenstein conjecture \cite{Fab99}, which predicts that $\sR(\sM_g)$ is a Gorenstein ring with socle in degree $g - 2$. The difficulty is to find enough relations between $\{\kappa_i\}$: starting from $g = 24$, all known methods have failed to do so. On the other hand, apart from numerical evidence, there has been no convincing reason for or against the conjecture.

One may consider variants of $\sR(\sM_g)$ for the universal curve $\sC_g$ and its relative powers~$\sC_g^n$. There is also a pointed version, \ie over the moduli of smooth pointed curves $\sM_{g, 1}$. In all cases one can ask about the corresponding Gorenstein properties.

\subsection{Jacobian side} --- Let $J$ be the Jacobian of a smooth pointed curve $(C, x_0)$. The point~$x_0$ induces an embedding $C \into J$. Since $J$ is an abelian variety, its Chow ring carries a second ring structure and also the action of the multiplication by $N \in \bZ$. The tautological ring $\sT(J)$ is defined to be the smallest $\bQ$-subspace of $\CH(J)$ (with $\bQ$-coeffcients) containing the curve class $[C]$, and stable under both ring structures as well as the multiplication by $N$. One can prove that $\sT(J)$ is finitely generated (with respect to the intersection product) and write down an explicit set of generators.

In a series of papers (\cite{Pol05}, \cite{Pol07}, \cite{Pol07b}, \etc), Polishchuk applied powerful tools to the study of $\sT(J)$, such as the Beauville decomposition, the Fourier transform, and the Lefschetz decomposition (or the $\fsl_2$-action). Notably he used the third tool to construct relations between the generators. For the generic curve (and modulo algebraic equivalence), he conjectured that they give all the relations (see \cite{Pol05}, Introduction).

Polishchuk's approach also brings a motivic touch to the subject. The hidden background is the so-called motivic Lefschetz decomposition, and the relations obtained this way are simply dictated by the motive of $J$.

\subsection{Main results} --- Our first result extends the Jacobian side to the relative (or universal) setting. Denote by $\sJ_{g, 1}$ the universal Jacobian over $\sM_{g, 1}$. The tautological ring~$\sT(\sJ_{g, 1})$ is defined similarly. By studying the $\fsl_2$-action on $\sT(\sJ_{g, 1})$, we obtain the following statement (Theorem~\ref{thm:structure}).

\secskip
\noindent\textit{Theorem 1.} --- \textit{The ring $\sT(\sJ_{g, 1})$ has an explicit finite set of generators \rleft with respect to the intersection product\rright. The $\fsl_2$-action on $\sT(\sJ_{g, 1})$ can also be described explicitly in terms of the generators.}
\secskip

Our second result connects the tautological rings on both sides. Denote by $\sC_{g, 1}$ the universal curve over $\sM_{g, 1}$, and by $\sC_{g, 1}^n$ (\resp $\sC_{g, 1}^{[n]}$) its $n$-th power (\resp symmetric power) relative to~$\sM_{g, 1}$. The tautological rings $\sR(\sC_{g, 1}^n)$ and $\sR(\sC_{g, 1}^{[n]})$ are defined accordingly. Further, one passes to the limit $\sC_{g, 1}^{[\infty]} = \varinjlim \sC_{g, 1}^{[n]}$ and get the tautological ring $\sR(\sC_{g, 1}^{[\infty]})$. Then we have a comparison (Corollary~\ref{cor:tautisoms} and Theorem~\ref{thm:poltaut}).

\secskip
\noindent\textit{Theorem 2.} --- \textit{The ring $\sR(\sC_{g, 1}^{[\infty]})$ is a polynomial algebra over $\sT(\sJ_{g, 1})$. In particular, the ring $\sR(\sM_{g, 1})$ is a $\bQ$-subalgebra of $\sT(\sJ_{g, 1})$.}
\secskip

Following Polishchuk, we use the $\fsl_2$-action on $\sT(\sJ_{g, 1})$ to produce tautological relations. With these relations, we are able to confirm that $\sR(\sM_{g, 1})$ (\resp $\sR(\sM_g)$) is Gorenstein for $g \leq 19$ (\resp $g \leq 23$). As far as computation goes, we seem to recover all the Faber-Zagier relations. Moreover, the socle condition for $\sR(\sC_{g, 1}^{[n]})$ allows us to formulate the corresponding Gorenstein property for $\sT(\sJ_{g, 1})$, which is again confirmed for $g \leq 7$. Our third result is an equivalence of Gorenstein properties (Theorem~\ref{thm:equiv}).

\secskip
\noindent\textit{Theorem 3.} --- \textit{The ring $\sT(\sJ_{g, 1})$ is Gorenstein if and only if $\sR(\sC_{g, 1}^{[n]})$ is Gorenstein for all $n \geq 0$.}
\secskip

It follows that $\sR(\sC_{g, 1}^{[n]})$ (\resp $\sR(\sC_g^{[n]})$) is Gorenstein for $g \leq 7$ and for all $n \geq 0$.

\subsectionempty We also spend a few words on the value of this note and how it relates to other works. By the work of Petersen and Tommasi (\cite{PT14} and \cite{Pet13}), the Gorenstein conjectures for the tautological rings of $\bar{\sM}_{g, n}$ and $\sM_{g, n}^{\ct}$ are known to be false already when $g = 2$. So it is tempting to believe that $\sR(\sM_g)$ may not be Gorenstein in general. Here we provide a possible explanation for this. Following Polishchuk's conjecture and philosophy in \cite{Pol05}, we may ask the following question.

\secskip
\noindent\textit{Question.} --- \textit{Are all relations of motivic nature?}
\secskip

This question, whose precise statement can be found in Conjecture~\ref{conj:motive}, has a certain geometric (rather than numerical) flavor. A positive answer to it would contradict the Gorenstein conjecture, and thus give an alternative description of the tautological rings.

More recently, a new set of relations for various tautological rings (including $\sR(\sC_g^{[n]})$) was conjectured by Pixton \cite{Pix12}, and proven by Pandharipande, Pixton and Zvonkine in cohomology \cite{PPZ13}, and by Janda in the Chow ring \cite{Jan13}. In \cite{Pix13}, Appendix~A, Pixton collected data about the discrepancies between his relations and the Gorenstein expectations for $\sR(\sC_g^{[n]})$. Our third result then shows that after $g = 24$ for $n = 0$, $g = 20$ for $n = 1$, \etc, the value $g = 8$ is the ultimate critical value for any large $n$. It is one of the most interesting cases to test various Gorenstein properties.

\subsection*{Conventions} --- We work over an arbitrary field $k$. Throughout, Chow rings $\CH = \oplus_i \CH^i$ are with $\bQ$-coefficients and graded by codimension. The symbol $\CH_i$ is only used relatively to a fixed base scheme $S$: if $X$ is a smooth projective scheme over $S$ with connected fibers, we write $\CH_i(X)$ for its Chow group of relative dimension~$i$ cycles with $\bQ$-coefficients. We set $\fsl_2 \coloneqq \bQ \cdot e + \bQ \cdot f + \bQ \cdot h$, with $[e, f] = h, [h, e] = 2e$ and $[h, f] = -2f$.

\subsection*{Acknowledgements} --- This note is part of the author's PhD thesis, written during his stays at the Université Paris VI, the University of Amsterdam and the Radboud University Nijmegen. It also expands an earlier preprint \cite{Yin12}.

The author is deeply indebted to his thesis advisor Ben Moonen, who introduced him to the subject and encouraged him throughout this project. He is  grateful to Li Ma, whose brilliant programming made this work much more meaningful. He also thanks Arnaud Beauville, Olof Bergvall, Alessandro Chiodo, Carel Faber, Gavril Farkas, Gerard van der Geer, Samuel Grushevsky, Felix Janda, Kefeng Liu, Rahul Pandharipande, Dan Petersen, Aaron Pixton, Sergey Shadrin, Mehdi Tavakol, Orsola Tommasi, Ravi Vakil, Claire Voisin (his co-advisor), Dmitry Zakharov, and Shengmao Zhu for many useful discussions and feedbacks. The whole project was inspired by Robin de Jong's talk at the PCMI in the summer of 2011.

The computation was carried out on the Lisa compute cluster at SURFsara (\href{http://surfsara.nl}{\texttt{surfsara.nl}}). The program code is available from the author upon request.

\secskip
\section{Tautological rings around a relative curve}
\secskip

\noindent We recall several tautological rings associated to a relative pointed curve. Then we focus on the universal situation, \ie over the moduli space, where (an analogue of) Faber's Gorenstein conjecture is stated. Throughout this section we work in the context of pointed curves, and we include a comparison with the unpointed counterpart at the end. 

\subsectionempty \label{sec:notation} Let $k$ be a field, and let $S$ be a smooth connected variety of dimension $d$ over $k$. Consider a relative curve $p \colon C \to S$ of genus $g$, \ie a smooth projective scheme over $S$ with geometrically connected fibers of relative dimension~$1$ and of genus~$g$. We assume $g > 0$. Further assume that $C/S$ admits a section (marked point) $x_0 \colon S \to C$.

For $n \geq 1$, denote by $p^n \colon C^n \to S$ (\resp $p^{[n]} \colon C^{[n]} \to S$) the $n$-th power (\resp symmetric power) of $C$ relative to $S$. Write $\sigma_n \colon C^n \to C^{[n]}$ for the symmetrization map. For convenience we set $C^0 = C^{[0]} \coloneqq S$.

\subsectionempty We describe a few geometric classes in the Chow rings $\CH(C^n)$, which serve as the building blocks of the tautological rings. First, denote by $K \in \CH^1(C)$ the first Chern class of the relative cotangent bundle $\Omega^1_{C/S}$. Also write $[x_0] \coloneqq \big[x_0(S)\big] \in \CH^1(C)$.

Next, define classes
\begin{align*}
\kappa_i & \coloneqq p_*(K^{i + 1}) \in \CH^i(S) \rfor i \geq 0, \\
\psi & \coloneqq x_0^*(K) \in \CH^1(S).
\end{align*}
We have $\kappa_0 = (2g - 2)[S]$, and it is convenient to write $\kappa_{-1} \coloneqq 0$. Also note that $x_0^*\big([x_0]\big) = -\psi$ by adjunction.

Further, we view $\CH(C^n)$ as $\CH(S)$-algebras by pulling back via $p^n$. We then keep the same notation $\{\kappa_i\}$ and $\psi$ for the pull-backs of these classes to $C^n$. For $1 \leq j \leq n$, let $\pr_j \colon C^n \to C$ be the projection to the $j$-th factor, and for $1 \leq k < l \leq n$, let $\pr_{k, l} \colon C^n \to C^2$ be the projection to the $k$-th and $l$-th factors. Denote by $\Delta$ the diagonal in $C^2$, and write
\begin{align*}
K_j & \coloneqq \pr_j^*(K) \in \CH^1(C^n), \\
[x_{0, j}] & \coloneqq \pr_j^*\big([x_0]\big) \in \CH^1(C^n), \\
[\Delta_{k, l}] & \coloneqq \pr_{k, l}^*\big([\Delta]\big) \in \CH^1(C^n).
\end{align*}

\subsection{Definition} --- Let $n \geq 0$. The \textit{tautological ring} of $C^n$, denoted by $\sR(C^n)$, is the $\bQ$-subalgebra of $\CH(C^n)$ generated by the classes $\{\kappa_i\}$, $\psi$, $\{K_j\}$ and $\big\{[x_{0, j}]\big\}$ (if $n \geq 1$), and $\big\{[\Delta_{k, l}]\big\}$ (if $n \geq 2$). Elements in $\sR(C^n)$ are called \textit{tautological classes}.

\secskip
In particular, the ring $\sR(S) \coloneqq \sR(C^0)$ is generated by $\{\kappa_i\}$ and $\psi$.

\subsection{Remark} \label{rem:deffp} --- Alternatively, one may define $\sR(C^n)$ in the style of Faber and Pandharipande (\cite{FP00} and \cite{FP05}). For $n \geq 1$ and $m \geq 0$, consider maps
\begin{equation*}
T = (T_1, \ldots, T_m) \colon C^n \to C^m,
\end{equation*}
such that each $T_i \colon C^n \to C$ is a projection of $C^n$ to one of its factors (when $m = 0$ we set $T = p^n \colon C^n \to S$). These maps are called \textit{tautological maps}. The tautological rings $\big\{\sR(C^n)\big\}$ then form the smallest system of $\bQ$-subalgebras satisfying $[x_0] \in \sR(C)$ and stable under pull-backs and push-forwards via all tautological maps. The proof is not difficult and is essentially in \cite{Loo95}, Proposition~2.1.

\subsectionempty For $n \geq 1$, the map $\sigma_n \colon C^n \to C^{[n]}$ induces an isomorphism of $\bQ$-algebras
\begin{equation*}
\sigma_n^* \colon \CH(C^{[n]}) \xto{\sim} \CH(C^n)^{\fS_n},
\end{equation*}
where $\CH(C^n)^{\fS_n}$ is the symmetric (or $\fS_n$-invariant) part of $\CH(C^n)$. Note that $\sigma_{n, *} \circ \sigma_n^* = n!$. The tautological ring of $C^{[n]}$ is then defined to be
\begin{equation*}
\sR(C^{[n]}) \coloneqq (\sigma_n^*)^{-1}\big(\sR(C^n)\big) = \frac{1}{n!}\sigma_{n, *}\big(\sR(C^n)\big).
\end{equation*}
When $n = 0$ we set $\sR(C^{[0]}) \coloneqq \sR(C^0)$ ($= \sR(S)$).

\subsection{Universal setting} --- The study of tautological rings was initiated by Mumford \cite{Mum83}, and later carried on extensively by Faber, Pandharipande, \etc, in the context of various moduli spaces of curves and compactifications. We refer to \cite{Fab99}, \cite{Pan02}, \cite{Fab13} and \cite{FP13} for an overview of major questions.

Our situation concerns the moduli space of smooth pointed curves of genus $g$ over $k$, denoted by $\sM_{g, 1}$ ($g > 0$ as before). It is isomorphic to the universal curve $\sC_g$ over the moduli of smooth genus $g$ curves $\sM_g$. We have $\dim(\sM_{g, 1}) = \dim(\sM_g) + 1 = 3g - 3 + 1$. The stacky nature of $\sM_{g, 1}$ does not play a role here. In fact, since $\sM_{g, 1}$ admits a finite cover by a smooth connected variety (see \cite{Mum83}, Part~I, Section~2), its Chow theory with $\bQ$-coefficients can be easily defined. In principle, one may regard $\sM_{g, 1}$ as a smooth connected variety in the sequel.

Denote by $\sC_{g, 1}$ the universal curve over $\sM_{g, 1}$. The tautological rings $\sR(\sC_{g, 1}^n)$ and $\sR(\sC_{g, 1}^{[n]})$ (and in particular $\sR(\sM_{g, 1})$) are thus defined by the same recipe as above, with $S = \sM_{g, 1}$.

\subsectionempty In \cite{Fab99}, Faber proposed the Gorenstein description of the tautological ring of $\sM_g$ (see Section~\ref{sec:unpoint} below). We adapt this idea to our context.

First, by the works of Looijenga \cite{Loo95} and Faber \cite{Fab97}, we have
\begin{equation} \label{eq:socle}
\sR^i(\sC_{g, 1}^n) = 0 \rfor i > g - 1 + n, \rand \sR^{g - 1 + n}(\sC_{g, 1}^n) \simeq \bQ.
\end{equation}
Then for $0 \leq i \leq g - 1 + n$, consider the pairing
\begin{equation} \label{eq:pairpower}
\sR^i(\sC_{g, 1}^n) \times \sR^{g - 1 + n - i}(\sC_{g, 1}^n) \xto{\cdot} \sR^{g - 1 + n}(\sC_{g, 1}^n) \simeq \bQ.
\end{equation}
The following is an analogue of Faber's Gorenstein conjecture.

\subsection{Speculation} \label{spec:faberpower} --- \textit{For $n \geq 0$ and $0 \leq i \leq g - 1 + n$, the pairing \textup{\eqref{eq:pairpower}} is perfect. In other words, the ring $\sR(\sC_{g, 1}^n)$ is Gorenstein with socle in degree $g - 1 + n$.}

\secskip
Note that $\sR^{g - 1 + n}(\sC_{g, 1}^n)$ is $\fS_n$-invariant. Restricting to the symmetric part, one would have the same Gorenstein property for $\sR(\sC_{g, 1}^{[n]})$, \ie there is a perfect pairing
\begin{equation} \label{eq:pairsym}
\sR^i(\sC_{g, 1}^{[n]}) \times \sR^{g - 1 + n - i}(\sC_{g, 1}^{[n]}) \xto{\cdot} \sR^{g - 1 + n}(\sC_{g, 1}^{[n]}) \simeq \bQ.
\end{equation}

\subsectionempty Speculation~\ref{spec:faberpower} is closely related to the Gorenstein conjecture for the tautological ring of $\sM_{g, n}^{\rt}$, the moduli space of stable $n$-pointed genus $g$ curves with rational tails. The latter is stated in \cite{Pan02}. In fact, for $n = 0$ we have exactly $\sC_{g, 1}^0 = \sM_{g, 1} = \sM_{g, 1}^{\rt}$. More generally there is a surjective map $\sM_{g, n + 1}^{\rt} \onto \sC_{g, 1}^n$, which can be expressed as a sequence of blow-ups. The Gorenstein properties on both sides are believed to be equivalent, but to the best of our knowledge there is not yet a written proof.

Moreover, Speculation~\ref{spec:faberpower} is proven for $g = 1, 2$ by Tavakol (\cite{Tav11} and \cite{Tav11b}; he went on proving the Gorenstein conjecture for the tautological rings of $\sM_{1, n}^{\rt}$ and $\sM_{2, n}^{\rt}$).

\subsection{Unpointed version} \label{sec:unpoint} --- Having stated everything in terms of pointed curves, we provide a translation to the classical, unpointed version.

Consider the universal curve $\sC_g$ over $\sM_g$ and its relative powers (\resp symmetric powers) $\sC_g^n$ (\resp $\sC_g^{[n]}$). The tautological rings $\sR(\sC_g^n)$ are defined to be generated by $\{\kappa_i\}$, $\{K_j\}$ and $\big\{[\Delta_{k, l}]\big\}$, \ie classes not involving the section $x_0$. One also defines $\sR(\sM_g)$ and $\sR(\sC_g^{[n]})$ accordingly.

Since $\sM_{g, 1} \simeq \sC_g$, there are isomorphisms
\begin{equation*}
\sC_{g, 1}^n \simeq \sC_g^n \times_{\sM_g} \sM_{g, 1} \simeq \sC_g^{n + 1}, \rand \sC_{g, 1}^{[n]} \simeq \sC_g^{[n]} \times_{\sM_g} \sM_{g, 1} \simeq \sC_g^{[n]} \times_{\sM_g} \sC_g.
\end{equation*}
Under these isomorphisms, we have a dictionary for tautological classes: the class $\psi \in \sR(\sC_{g, 1}^n)$ corresponds to $K_{n + 1} \in \sR(\sC_g^{n + 1})$, and $[x_{0, j}] \in \sR(\sC_{g, 1}^n)$ corresponds to $[\Delta_{j, n + 1}] \in \sR(\sC_g^{n + 1})$. The dictionary gives isomorphisms of $\bQ$-algebras
\begin{equation*}
\sR(\sC_{g, 1}^n) \simeq \sR(\sC_g^{n + 1}), \rand \sR(\sC_{g, 1}^{[n]}) \simeq \sR(\sC_g^{n + 1})^{\fS_n},
\end{equation*}
where $\fS_n$ acts on the first $n$ factors of $\sC_g^{n + 1}$. 

One can also formulate the corresponding Gorenstein properties for $\sR(\sC_g^n)$ and $\sR(\sC_g^{[n]})$ (note that $\sR^{g - 2 + n}(\sC_g^n) = \sR^{g - 2 + n}(\sC_g^{[n]}) \simeq \bQ$). It is immediate that the Gorenstein properties for $\sR(\sC_{g, 1}^n)$ and $\sR(\sC_g^{n + 1})$ are equivalent, and that $\sR(\sC_{g, 1}^{[n]})$ being Gorenstein implies $\sR(\sC_g^{[n + 1]})$ being Gorenstein.

\secskip
\section{The Chow ring of the relative Jacobian}
\secskip

\noindent We review three important structures on the Chow ring of the relative Jacobian, namely the \textit{Beauville decomposition}, the \textit{Fourier transform} and the \textit{Lefschetz decomposition}. We also discuss their motivic background. A picture called the Dutch house is presented to visualize these structures.

\subsectionempty \label{sec:notation2} As in the previous section $S$ is a smooth connected variety of dimension $d$ over $k$. Let $p \colon C \to S$ be a relative curve of genus $g > 0$, together with a section $x_0 \colon S \to C$. Denote by $\pi \colon J \to S$ the associated relative Jacobian. By definition $J \coloneqq \Pic^0(C/S)$ is an abelian scheme of relative dimension~$g$. The section $x_0$ induces an embedding $\iota \colon C \into J$, which sends (locally) a section $x$ of $C/S$ to $\iota(x) \coloneqq \sO_C(x - x_0)$. The composition $\iota \circ x_0$ is then the zero section $o \colon S \to J$. To summarize, we have the following diagram.
\begin{equation} \label{eq:setting}
\begin{tikzcd}[column sep=small]
C \arrow[hookrightarrow]{rr}{\iota} \arrow{dr}[pos=.56]{p} & & J \arrow{dl}[swap]{\pi} \\
& S \arrow[bend left=50]{ul}[pos=.51]{x_0} \arrow[bend right=50]{ur}[swap]{o}
\end{tikzcd}
\end{equation}

\subsectionempty The (abelian) group structure on $J$ gives the addition and scalar multiplication maps
\begin{equation*}
\mu \colon J \times_S J \to J, \rand [N] \colon J \to J \rfor N \in \bZ.
\end{equation*}
The Chow ring $\CH(J)$ carries a second ring structure called the \textit{Pontryagin product} ``$*$''. It sends $\alpha \in \CH^i(J)$ and $\beta \in \CH^j(J)$ to
\begin{equation} \label{eq:pontryagin}
\alpha * \beta \coloneqq \mu_*\big(\pr_1^*(\alpha) \cdot \pr_2^*(\beta)\big) \in \CH^{i + j - g}(J),
\end{equation}
where $\pr_1, \pr_2 \colon J \times_S J \to J$ are the two projections.

\subsection{Beauville decomposition} --- By the work of Beauville \cite{Bea86} (later generalized to the relative setting by Deninger and Murre \cite{DM91}), the Chow ring $\CH(J)$ can be decomposed into eigenspaces with respect to the action of $[N]$. More precisely, for $0 \leq i \leq g + d$ we have
\begin{equation} \label{eq:beauville}
\CH^i(J) = \bigoplus_{j = \max\{i - g, 2i - 2g\}}^{\min\{i + d, 2i\}} \CH^i_{(j)}(J),
\end{equation}
with
\begin{equation*}
\CH^i_{(j)}(J) \coloneqq \big\{\alpha \in \CH^i(J) : [N]^*(\alpha) = N^{2i - j} \alpha \rforall N \in \bZ\big\}.
\end{equation*}
The decomposition is multiplicative: for $\alpha \in \CH^i_{(j)}(J)$ and $\beta \in \CH^k_{(l)}(J)$ we have $\alpha \cdot \beta \in \CH^{i + k}_{(j + l)}(J)$ and $\alpha * \beta \in \CH^{i + k - g}_{(j + l)}(J)$.

\subsectionempty To introduce the two other structures, recall that $J$ has a canonical principal polarization $\lambda \colon J \xto{\sim} J^t$. Denote by $\sP$ (the pull-back of) the Poincaré line bundle on $J \times_S J$, trivialized along the two zero sections. We write $\ell \in \CH^1(J \times_S J)$ for its first Chern class.

Also associated to the polarization is a relatively ample divisor class $\theta \in \CH^1_{(0)}(J)$. There are identities
\begin{equation} \label{eq:poincare}
\theta = -\frac{1}{2}\Delta^*(\ell), \rand \ell = \pr_1^*(\theta) + \pr_2^*(\theta) - \mu^*(\theta),
\end{equation}
where $\Delta \colon J \to J \times_S J$ is the diagonal map. Note that we adopted the sign convention explained in \cite{Bea10}, Section~1.6.

\subsection{Fourier transform} --- The class $\exp(\ell)$, viewed as a correspondence, induces the Fourier transform $\sF \colon \CH(J) \to \CH(J)$. It sends $\alpha \in \CH(J)$ to
\begin{equation*}
\sF(\alpha) \coloneqq \pr_{2, *}\big(\pr_1^*(\alpha) \cdot \exp(\ell)\big).
\end{equation*}
We have the following identities (see \cite{Bea86}, Section~1 or \cite{DM91}, Section~2)
\begin{equation*}
\sF(\alpha * \beta) = \sF(\alpha) \cdot \sF(\beta), \ \ \sF(\alpha \cdot \beta) = (-1)^g\sF(\alpha) * \sF(\beta), \rand \sF \circ \sF = (-1)^g[-1]^*.
\end{equation*}
Further, there are isomorphisms
\begin{equation*}
\sF \colon \CH^i_{(j)}(J) \xto{\sim} \CH^{g - i + j}_{(j)}(J).
\end{equation*}
Hence by applying $\sF$ to $\CH^i(J)$ and by collecting components of different codimensions, one recovers the Beauville decomposition.

\subsection{Lefschetz decomposition} --- By Künnemann \cite{Kun93}, the classical Lefschetz decomposition in cohomology can be lifted to the Chow ring of an abelian scheme. In the case of the relative Jacobian, Polishchuk \cite{Pol07b} showed that the decomposition can be reconstructed geometrically.

Write $[C] \coloneqq \big[\iota(C)\big] \in \CH^{g - 1}(J)$ for the curve class, and $[C]_{(j)} \in \CH^{g - 1}_{(j)}(J)$ for the components of $[C]$ in the Beauville decomposition. First, we have identities
\begin{equation} \label{eq:deftheta}
\theta = \frac{1}{(g - 1)!}[C]_{(0)}^{*(g - 1)} = -\sF\big([C]_{(0)}\big).
\end{equation}
On $\CH(J)$, define operators
\begin{align}
e \colon \CH^i_{(j)}(J) & \to \CH^{i + 1}_{(j)}(J) & \alpha & \mapsto -\theta \cdot \alpha, \nonumber \\
f \colon \CH^i_{(j)}(J) & \to \CH^{i - 1}_{(j)}(J) & \alpha & \mapsto -[C]_{(0)} * \alpha, \label{eq:sl2}\\
h \colon \CH^i_{(j)}(J) & \to \CH^i_{(j)}(J) & \alpha & \mapsto (2i - j - g)\alpha. \nonumber
\end{align}
Then the operators $e, f$ and $h$ generate a $\bQ$-linear representation of $\fsl_2$, which induces the Lefschetz decomposition on $\CH(J)$. See \cite{Pol07b}, Theorem~2.6 for the proof of this statement. Note that we followed Polishchuk's sign convention by setting $e$ to be $-\theta$.

From now on we refer to \eqref{eq:sl2} as \textit{the $\fsl_2$-action}.

\subsection{Remark} \label{rem:exp} --- There is an identity (see \cite{Bea04}, Section~2.3~(iv) or \cite{Pol08}, Lemma~1.4)
\begin{equation*}
\sF = \exp(e)\exp(-f)\exp(e) \ron \CH(J).
\end{equation*}
This means if we represent $e, f$ by the matrices $\big(\begin{smallmatrix} 0 & 1 \\ 0 & 0 \end{smallmatrix}\big), \big(\begin{smallmatrix} 0 & 0 \\ 1 & 0 \end{smallmatrix}\big) \in \fsl_2(\bQ)$, then $\sF$ corresponds to the matrix $\big(\begin{smallmatrix} 0 & 1 \\ -1 & 0 \end{smallmatrix}\big) \in \SL_2(\bQ)$. As the Fourier transform induces the Beauville decomposition, we see that among the three structures above, the Lefschetz decomposition (or the $\fsl_2$-action) may be viewed as the deepest.

\subsection{Motivic background} --- We briefly discuss the motivic interpretation of the structures above. See \cite{DM91} and \cite{Kun93} for more details.

Denote by $R(J/S)$ the relative Chow motive of $J$. Then $R(J/S)$ admits a decomposition
\begin{equation*}
R(J/S) = \bigoplus_{i = 0}^{2g}R^i(J/S),
\end{equation*}
with $[N]^*$ acting on $R^i(J/S)$ by multiplication by $N^i$. It follows that
\begin{equation*}
\CH^i_{(j)}(J) = \CH^i\big(R^{2i - j}(J/S)\big).
\end{equation*}
Moreover, the Fourier transform $\sF$ induces isomorphisms
\begin{equation*}
\sF \colon R^i(J/S) \xto{\sim} R^{2g - i}(J/S)(g - i),
\end{equation*}
where ``$(\textup{-})$'' stands for Tate twists. Finally, we have the \textit{motivic Lefschetz decomposition}
\begin{equation*}
R^i(J/S) = \bigoplus_{k = \max\{0, i - g\}}^{\lfloor i/2 \rfloor} e^k\big(R^{i - 2k}_{\prim}(J/S)(-k)\big),
\end{equation*}
where $R^i_{\prim}(J/S)$ is the primitive part of $R^i(J/S)$ with respect to the $\fsl_2$-action.

\subsectionempty We find it convenient to replace the codimension grading on $\CH(J)$ by a new, motivic grading. We write
\begin{equation*}
\CH_{(i, j)}(J) \coloneqq \CH^{(i + j)/2}_{(j)}(J), \textup{ \ or equivalently \ } \CH_{(2i - j, j)}(J) \coloneqq \CH^i_{(j)}(J),
\end{equation*}
so that $[N]^*$ acts on $\CH_{(i, j)}(J)$ by multiplication by $N^i$. In other words, we set
\begin{equation*}
\CH_{(i, j)}(J) = \CH\big(R^i(J/S)\big) \cap \CH_{(j)}(J).
\end{equation*}
The Beauville decomposition \eqref{eq:beauville} then takes the form
\begin{equation} \label{eq:beauville2}
\CH(J) = \bigoplus_{i, j} \CH_{(i, j)}(J),
\end{equation}
with $0 \leq i \leq 2g$, $\max\{-i, i - 2g\} \leq j \leq \min\{i, 2g - i\} + 2d$ and $i + j$ even.

Expressions under the new grading are simple (if not simpler): for $\alpha \in \CH_{(i, j)}(J)$ and $\beta \in \CH_{(k, l)}(J)$ we have $\alpha \cdot \beta \in \CH_{(i + k, j + l)}(J)$ and $\alpha * \beta \in \CH_{(i + k - 2g, j + l)}(J)$. Also $\sF$ induces
\begin{equation*}
\sF \colon \CH_{(i, j)}(J) \xto{\sim} \CH_{(2g - i, j)}(J).
\end{equation*}
Further, we have $\theta \in \CH_{(2, 0)}(J)$. The $\fsl_2$-action in \eqref{eq:sl2} becomes
\begin{align*}
e \colon \CH_{(i, j)}(J) & \to \CH_{(i + 2, j)}(J) & \alpha & \mapsto -\theta \cdot \alpha, \\
f \colon \CH_{(i, j)}(J) & \to \CH_{(i - 2, j)}(J) & \alpha & \mapsto -[C]_{(0)} * \alpha, \\
h \colon \CH_{(i, j)}(J) & \to \CH_{(i, j)}(J) & \alpha & \mapsto (i - g)\alpha.
\end{align*}

\subsection{Dutch house} --- We present a useful picture that illustrates the structures above and combines the motivic aspect. It also allows us to make clear statements without complicated indices. The picture is inspired by \cite{Moo09}, Figure~1. We call it the \textit{Dutch house} due to its resemblance to a traditional Dutch \textit{trapgevel}.

\begin{figure}
\centering
\includegraphics[height=.7\textheight]{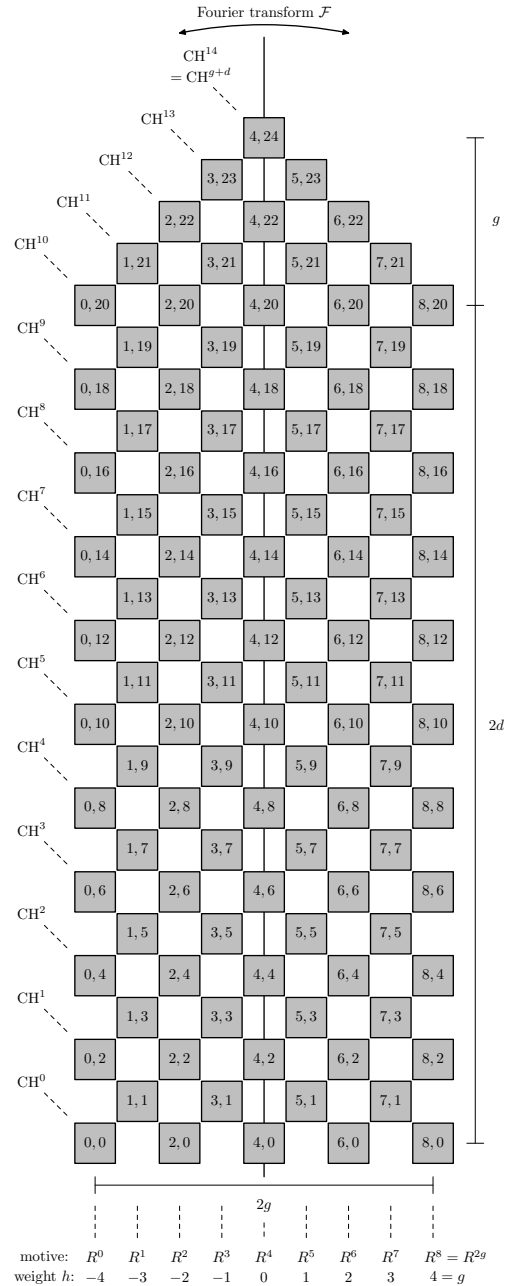}
\caption{Mon dessin n\textsuperscript{o} 1: the outside of the Dutch house ($g = 4$ and $d = 10$).}
\label{fig:outside}
\end{figure}

In Figure~\ref{fig:outside}, the $(i, j)$-th block represents the component $\CH_{(i, j)}(J)$ in the Beauville decomposition. The columns then read the motivic decomposition $R(J/S) = \oplus_i R^i(J/S)$, and the rows read Beauville's grading $j$. As a result, components with the same codimension lie on a dashed line from upper left to lower right.

It is not difficult to verify that the house shape results from the precise index range of~\eqref{eq:beauville2}. The width of the house depends on the genus~$g$, while the height (without roof) depends on~$d = \dim(S)$. In particular when $S = k$, the house reduces to the roof only. Here Figure~\ref{fig:outside} is drawn based on the universal Jacobian over $S = \sM_{4, 1}$, with $d = \dim(\sM_{4, 1}) = 10$.

As is shown in the picture, the Fourier transform $\sF$ acts as the reflection over the middle vertical line. The operator $e$ in the $\fsl_2$-action shifts classes to the right by $2$ blocks, while $f$ shifts classes to the left by $2$ blocks. Finally the middle column of the house has weight $0$ with respect to $h$.

\subsection{Remark} --- Note that we have not drawn the components $\CH_{(i, j)}(J)$ with negative $j$. On one hand, when $S = k$ the Beauville conjecture predicts the vanishing of those components (see \cite{Bea86}, Section~2). On the other hand, all classes we shall encounter are in $\CH_{(i, j)}(J)$ with $j \geq 0$, \ie inside the house.

\subsectionempty The isomorphisms of motives (see \cite{DM91}, Example~1.4)
\begin{equation*}
\begin{tikzcd}[column sep=tiny,row sep=scriptsize]
\phantom{(g^g)} R^0(J/S) \arrow[leftarrow]{rr}{\sF}[swap]{\sim} & & R^{2g}(J/S)(g) \arrow{dl}[pos=.61]{\ \pi_*}[swap]{\sim} \\
& R(S/S) \arrow{ul}{\pi^*}[swap]{\sim}
\end{tikzcd}
\end{equation*}
induce the following isomorphisms of $\bQ$-algebras.
\begin{equation} \label{eq:chowisoms}
\begin{tikzcd}[column sep=tiny,row sep=scriptsize]
\phantom{_g} \big(\bigoplus_{i = 0}^{d}\CH_{(0, 2i)}(J), \cdot\big) \arrow[leftarrow]{rr}{\sF}[swap]{\sim} & & \big(\bigoplus_{i = 0}^{d}\CH_{(2g, 2i)}(J), *\big) \arrow{dl}[pos=.61]{\ \pi_*}[swap]{\sim} \\
& \big(\CH(S), \cdot\big) \arrow{ul}{\pi^*}[swap]{\sim}
\end{tikzcd}
\end{equation}
The gradings are preserved as $\pi^* \colon \CH^i(S) \xto{\sim} \CH_{(0, 2i)}(J)$ and $\pi_* \colon \CH_{(2g, 2i)}(J) \xto{\sim} \CH^i(S)$.

In particular, the Chow ring $\CH(S)$ may be regarded as a $\bQ$-subalgebra of $\big(\CH(J), \cdot\big)$ via~$\pi^*$, or as a $\bQ$-subalgebra of $\big(\CH(J), *\big)$ via~$\pi_*$. In terms of the Dutch house, we may identify $\CH(S)$ with the $0$-th column or with the $2g$-th column of the house.

\secskip
\section{The tautological ring of the relative Jacobian}
\secskip

\noindent We define the tautological ring $\sT(J)$ of the relative Jacobian. We determine its generators and give explicit formulae for the $\fsl_2$-action in terms of the generators. The consequence is that by pulling back via $\pi^* \colon \CH(S) \to \CH(J)$, one can identify $\sR(S)$ with the $\bQ$-subalgebra of~$\sT(J)$ located on the $0$-th column of the Dutch house. Throughout this section we work in the setting of \eqref{eq:setting}.

\subsection{Definition} --- The \textit{tautological ring} of $J$, denoted by $\sT(J)$, is the smallest $\bQ$-subalgebra of $\big(\CH(J), \cdot, *\big)$ (\ie subalgebra with respect to both ``$\cdot$'' and ``$*$'') that contains $[C] \in \CH^{g - 1}(J)$ and that is stable under $[N]^*$, for all $N \in \bZ$. Elements in $\sT(J)$ are called \textit{tautological classes}.

\secskip
This notion of tautological ring was introduced by Beauville \cite{Bea04} in the context of a Jacobian variety and modulo algebraic equivalence. Since then it has been studied in various contexts. We refer to \cite{Pol05}, \cite{Pol07}, \cite{Her07}, \cite{GK07} and \cite{Moo09} for more details. In the relative setting, Polishchuk considered a much bigger tautological ring, including all classes in~$\pi^*\big(\CH(S)\big)$ (see \cite{Pol07b}, Section~4). Our minimalist version turns out to be the right one for studying the tautological ring of $S$.

\subsection{Remarks} --- (i) \ The ring $\sT(J)$ is stable under the Beauville decomposition, the Fourier transform and the Lefschetz decomposition (or the $\fsl_2$-action). Indeed, by definition $\sT(J)$ is graded by codimension. Also by applying $[N]^*$ with various $N$ we have $[C]_{(0)} \in \sT(J)$. Then it follows from~\eqref{eq:deftheta} and \eqref{eq:sl2} that $\sT(J)$ is stable under the $\fsl_2$-action. As is shown in Remark~\ref{rem:exp}, this implies the stability of $\sT(J)$ under all three structures on $\CH(J)$.

(ii) \ Alternatively, one may define $\sT(J)$ to be the smallest $\bQ$-subalgebra of $\big(\CH(J), \cdot\big)$ containing $[C] \in \CH^{g - 1}(J)$ and stable under the $\fsl_2$-action (and thus stable under the other two structures). In fact, being stable under the Beauville decomposition is the same as being stable under $[N]^*$. Since $\sF$ interchanges ``$\cdot$'' and ``$*$'', the ring is also stable under ``$*$''.

\subsectionempty Since the two products ``$\cdot$'' and ``$*$'' do not commute with each other, it is \apriori not clear whether $\sT(J)$ is finitely generated. We give an affirmative answer to this question by writing down an explicit set of generators.

Recall that $[C]_{(j)} \in \sT_{(2g - 2 - j, j)}(J)$. By the index range of \eqref{eq:beauville2}, we have $[C]_{(j)} = 0$ for $j < 0$ or $j > 2g - 2$. Now consider for $i \leq j + 2$ and $i + j$ even
\begin{equation*}
\theta^{(j - i + 2)/2} \cdot [C]_{(j)} \in \sT_{(2g - i, j)}(J).
\end{equation*}
Denote its Fourier dual by
\begin{equation*}
p_{i, j} \coloneqq \sF\big(\theta^{(j - i + 2)/2} \cdot [C]_{(j)}\big) \in \sT_{(i, j)}(J).
\end{equation*}
As examples we have $p_{2, 0} = \sF\big([C]_{(0)}\big) = -\theta$ and $p_{0, 0} = \sF\big(\theta \cdot [C]_{(0)}\big) = g[J]$. The index range of \eqref{eq:beauville2} implies $p_{i, j} = 0$ for $i < 0$ or $j < 0$ or $j > 2g - 2$.

\begin{figure}
\centering
\includegraphics[height=.5\textheight]{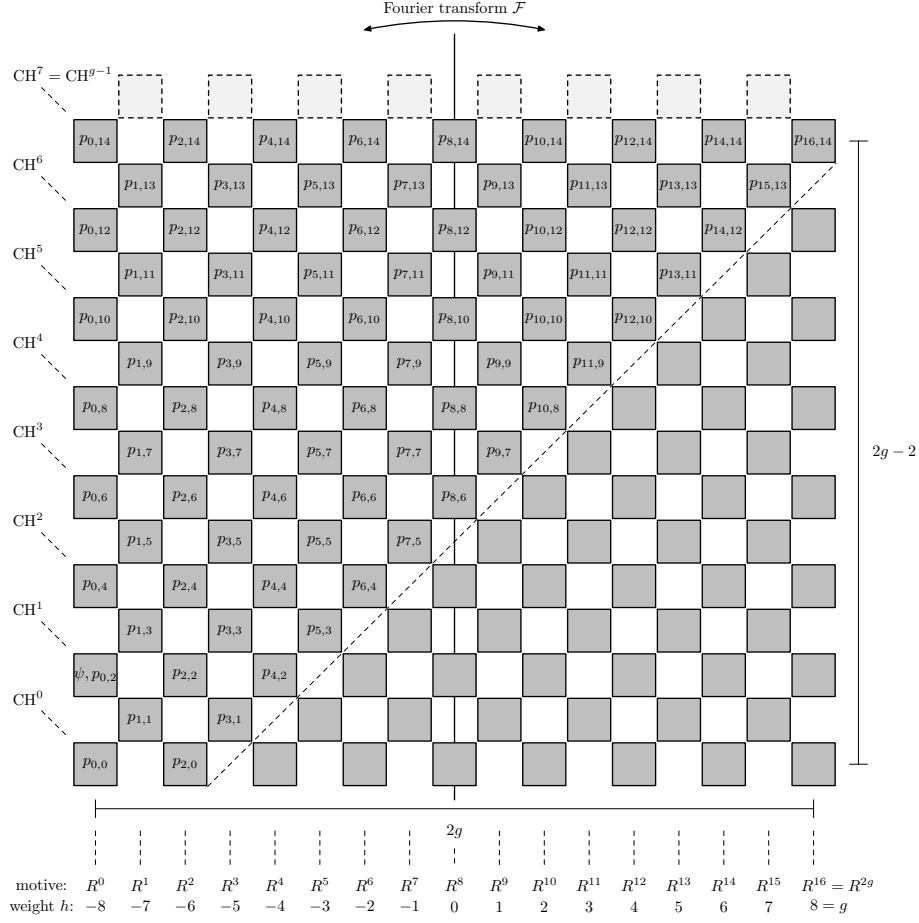}
\caption{Mon dessin n\textsuperscript{o} 2: the inside of the Dutch house ($g = 8$).}
\label{fig:inside}
\end{figure}

Figure~\ref{fig:inside} depicts the classes $\{p_{i, j}\}$ inside the Dutch house with $g = 8$. Also shown in the picture is the pull-back of the class $\psi$ via $\pi^*$, again denoted by $\psi$, which lies in $\CH_{(0, 2)}(J)$. Note that when $d = \dim(S)$ is small, classes above the roof vanish as well.

\subsectionempty By definition, the operator $e$ in the $\fsl_2$-action is the intersection with $p_{2, 0}$. Also it is not difficult to verify that
\begin{equation*}
f(p_{i, j}) = p_{i - 2, j}.
\end{equation*}
One of the questions is to calculate the class $f(p_{i, j} p_{k, l})$. This turns out to be the key to the following theorem.

\subsection{Theorem} \label{thm:structure} --- (i) \ \textit{The ring $\sT(J)$ coincides with the $\bQ$-subalgebra of $\big(\CH(J), \cdot\big)$ generated by the classes $\{p_{i, j}\}$ and $\psi$. In particular, it is finitely generated.}

(ii) \ \textit{The operator $f$ acts on polynomials in $\{p_{i, j}\}$ and $\psi$ via the degree $2$ differential operator}
\begin{equation} \label{eq:defd}
\begin{split}
\sD \coloneqq \frac{1}{2}\sum_{i, j, k, l}\Bigg(\psi p_{i - 1, j - 1} p_{k - 1, l - 1} & - \binom{i + k - 2}{i - 1} p_{i + k - 2, j + l}\Bigg) \partial p_{i, j} \partial p_{k, l} \\
& + \sum_{i, j}p_{i - 2, j}\partial p_{i, j}. 
\end{split}
\end{equation}

\secskip
To prove the theorem, we begin with a lemma that computes the pull-back of $\theta$ via the embedding $\iota \colon C \into J$. It is probably known to experts, and is implicit in \cite{Pol07b}, Theorem~2.6.

\subsection{Lemma} \label{lem:pulltheta} --- \textit{We have an identity}
\begin{equation} \label{eq:pulltheta}
\iota^*(\theta) = \frac{1}{2}K + [x_0] + \frac{1}{2}\psi \rin \CH^1(C).
\end{equation}

\begin{proof}
The goal is to calculate $\iota^*(\theta) = -\iota^*\Big(\sF\big([C]_{(0)}\big)\Big)$ and we start from $\iota^*\Big(\sF\big([C]\big)\Big)$. Consider the following three Cartesian squares
\begin{equation*}
\begin{tikzcd}
C \times_S C \arrow{r}{\iota \times_S \id_C} \arrow{d}[swap]{\id_C \times_S \iota} & J \times_S C \arrow{r}{\pr_2} \arrow{d}{\id_J \times_S \iota} & C \arrow{d}{\iota} \\
C \times_S J \arrow{r}[swap]{\iota \times_S \id_J} \arrow{d}[swap]{\pr_1} & J \times_S J \arrow{r}[swap]{\pr_2} \arrow{d}{\pr_1} & J \\
C \arrow{r}[swap]{\iota} & J & \phantom{C \times_S C}
\end{tikzcd}
\end{equation*}
where $\pr_1$ and $\pr_2$ stand for the two projections in all cases. Then we have
\begin{align*}
\iota^*\Big(\sF\big([C]\big)\Big) & = \iota^*\pr_{2, *}\Big(\pr_1^*\,\iota_*\big([C]\big) \cdot \exp(\ell)\Big) \\
& = \pr_{2, *}(\id_J \times_S \iota)^*\Big((\iota \times_S \id_J)_*\pr_1^*\big([C]\big) \cdot \exp(\ell)\Big) \displaybreak\\
& = \pr_{2, *}(\id_J \times_S \iota)^*\Big((\iota \times_S \id_J)_*\big([C \times_S J]\big) \cdot \exp(\ell)\Big) \\
& = \pr_{2, *}(\id_J \times_S \iota)^*(\iota \times_S \id_J)_*(\iota \times_S \id_J)^*\big(\exp(\ell)\big) \\
& = \pr_{2, *}(\iota \times_S \id_C)_*(\id_C \times_S \iota)^*(\iota \times_S \id_J)^*\big(\exp(\ell)\big) \\
& = \pr_{2, *}(\iota \times_S \iota)^*\big(\exp(\ell)\big) \\
& = \pr_{2, *}\Big(\exp\big((\iota \times_S \iota)^*(\ell)\big)\Big).
\end{align*}

The second identity in \eqref{eq:poincare} and the theorem of the square imply (see \cite{Pol07b}, Formula~(2.1))
\begin{equation} \label{eq:pulll}
(\iota \times_S \iota)^*(\ell) = [\Delta] - \pr_1^*\big([x_0]\big) - \pr_2^*\big([x_0]\big) - \psi,
\end{equation}
where $\Delta$ is the diagonal in $C \times_S C$. It follows that
\begin{equation} \label{eq:pullfc}
\begin{split}
\iota^*\Big(\sF\big([C]\big)\Big) & = \pr_{2, *}\bigg(\exp\Big([\Delta] - \pr_1^*\big([x_0]\big) - \pr_2^*\big([x_0]\big) - \psi\Big)\bigg) \\
& = \pr_{2, *}\bigg(\exp\Big([\Delta] - \pr_1^*\big([x_0]\big)\Big)\bigg) \cdot \exp\big({-[x_0]} - \psi\big).
\end{split}
\end{equation}
Observe that on the left-hand side of (\ref{eq:pullfc}), we have
\begin{equation*}
\iota^*\Big(\sF\big([C]\big)\Big) = \sum_{j = 0}^{2g - 2} \iota^*\Big(\sF\big([C]_{(j)}\big)\Big),
\end{equation*}
with $\iota^*\Big(\sF\big([C]_{(j)}\big)\Big) \in \CH^{j + 1}(C)$. Hence $\iota^*\Big(\sF\big([C]_{(0)}\big)\Big)$ is just the codimension~$1$ component of $\iota^*\Big(\sF\big([C]\big)\Big)$. Expanding the exponentials in (\ref{eq:pullfc}) while keeping track of the codimension, we obtain
\begin{align*}
\iota^*\Big(\sF\big([C]_{(0)}\big)\Big) & = \pr_{2, *}\bigg(\frac{1}{2}\Big([\Delta] - \pr_1^*\big([x_0]\big)\Big)^2\bigg) - \pr_{2, *}\Big([\Delta] - \pr_1^*\big([x_0]\big)\Big) \cdot \big([x_0] + \psi\big) \\
& = \pr_{2, *}\bigg(\frac{1}{2}\Big([\Delta] - \pr_1^*\big([x_0]\big)\Big)^2\bigg) \\
& = \frac{1}{2}\pr_{2, *}\big([\Delta] \cdot [\Delta]\big) - \pr_{2, *}\Big([\Delta] \cdot \pr_1^*\big([x_0]\big)\Big) + \frac{1}{2}\pr_{2, *}\pr_1^*\big([x_0] \cdot [x_0]\big).
\end{align*}

The first two terms in the previous expression are easily calculated: we have
\begin{equation*}
\pr_{2, *}\big([\Delta] \cdot [\Delta]\big) = -K, \rand \pr_{2, *}\Big([\Delta] \cdot \pr_1^*\big([x_0]\big)\Big) = [x_0].
\end{equation*}
For the term $\pr_{2, *}\pr_1^*\big([x_0] \cdot [x_0]\big)$, consider the following Cartesian square.
\begin{equation*}
\begin{tikzcd}
C \times_S C \arrow{r}{\pr_2} \arrow{d}[swap]{\pr_1} & C \arrow{d}{p} \\
C \arrow{r}[swap]{p} & S
\end{tikzcd}
\end{equation*}
Then we have
\begin{equation*}
\pr_{2, *}\pr_1^*\big([x_0] \cdot [x_0]\big) = p^*p_*\big([x_0] \cdot [x_0]\big) = p^*p_*x_{0, *}x_0^*\big([x_0]\big) = p^*x_0^*\big([x_0]\big) = -\psi.
\end{equation*}
In total we find $\iota^*\Big(\sF\big([C]_{(0)}\big)\Big) = -K/2 - [x_0] - \psi/2$.
\end{proof}

\subsection{Proof of Theorem~\textup{\ref{thm:structure}}} --- Suppose we have proven (ii) and that $\psi \in \sT(J)$. Consider the $\bQ$-subalgebra of $\big(\CH(J), \cdot\big)$ generated by the classes $\{p_{i, j}\}$ and $\psi$. We denote it by $\sT'(J)$ and we have $\sT'(J) \subset \sT(J)$. By definition $\sT'(J)$ is stable under the action of $e \in \fsl_2$. Now (ii) shows that $\sT'(J)$ is also stable under the action of $f \in \fsl_2$. It then follows from Remark~\ref{rem:exp} that $\sT'(J)$ is stable under $\sF$. In particular, the classes $\big\{[C]_{(j)}\big\}$ are in $\sT'(J)$. Since $\sT(J)$ is the smallest $\bQ$-subalgebra containing $[C]$ and stable under the $\fsl_2$-action, there is necessarily an equality $\sT'(J) = \sT(J)$, which proves (i).

Statement~(ii) follows essentially from \cite{Pol07b}, Formula~(2.9). We just need to translate the notation carefully. Following Polishchuk, we write $\eta \coloneqq K/2 + [x_0] + \psi/2$, which by \eqref{eq:pulltheta} is equal to $\iota^*(\theta)$. We also have $f = -\tilde{X}_{2, 0}(C)/2$ in his notation. Define operators $\tilde{p}_{i, j}$ on $\CH(J)$ by $\tilde{p}_{i, j}(\alpha) \coloneqq p_{i, j} \cdot \alpha$. Then the fact that 
\begin{equation*}
p_{i, j} = \sF\big(\theta^{(j - i + 2)/2} \cdot [C]_{(j)}\big) = \sF\big(\iota_*(\eta^{(j - i + 2)/2})_{(j)}\big)
\end{equation*}
is translated into
\begin{equation*}
\tilde{p}_{i, j} = \frac{1}{i!}\tilde{X}_{0, i}(\eta^{(j - i + 2)/2}).
\end{equation*}
We apply Formula~(2.9) in \loccit~and find
\begin{align*}
[f, \tilde{p}_{i, j}] & = -\frac{1}{2 \cdot i!}\big[\tilde{X}_{2, 0}(C), \tilde{X}_{0, i}(\eta^{(j - i + 2)/2})\big] \\
& = \frac{1}{(i - 1)!}\tilde{X}_{1, i - 1}(\eta^{(j - i + 2)/2}) - \frac{1}{(i - 2)!}\tilde{X}_{0, i - 2}(\eta^{(j - i + 4)/2}).
\end{align*}
Note that the second equality above also involves the fact that $\tilde{X}_{i, 0}(C) = 0$ for $i \leq 1$ (see \loccit, Lemma~2.8), and that $x_0^*(\eta) = x_0^*\iota^*(\theta) = o^*(\theta) = 0$. We continue to calculate
\begin{align*}
\big[[f, \tilde{p}_{i, j}], \tilde{p}_{k, l}\big] & = \frac{1}{(i - 1)!k!}\big[\tilde{X}_{1, i - 1}(\eta^{(j - i + 2)/2}), \tilde{X}_{0, k}(\eta^{(l - k + 2)/2})\big] \\
& \qquad - \frac{1}{(i - 2)!k!}\big[\tilde{X}_{0, i - 2}(\eta^{(j - i + 4)/2}), \tilde{X}_{0, k}(\eta^{(l - k + 2)/2})\big].
\end{align*}
Applying the same formula, we have $\big[\tilde{X}_{0, i - 2}(\eta^{(j - i + 4)/2}), \tilde{X}_{0, k}(\eta^{(l - k + 2)/2})\big] = 0$, and
\begin{align*}
\big[\tilde{X}_{1, i - 1}(\eta^{(j - i + 2)/2}), \tilde{X}_{0, k}(\eta^{(l - k + 2)/2})\big] & = k \psi \tilde{X}_{0, k - 1}(\eta^{(l - k + 2)/2}) \tilde{X}_{0, i - 1}(\eta^{(j - i + 2)/2}) \\
& \qquad -k \tilde{X}_{0, i + k - 2}(\eta^{(j - i + l - k + 4)/2}).
\end{align*}
In total, we obtain
\begin{equation} \label{eq:deg2}
\begin{split}
\big[[f, \tilde{p}_{i, j}], \tilde{p}_{k, l}\big] & = \frac{1}{(i - 1)!(k - 1)!}\Big(\psi \tilde{X}_{0, k - 1}(\eta^{(l - k + 2)/2}) \tilde{X}_{0, i - 1}(\eta^{(j - i + 2)/2}) \\
& \qquad - \tilde{X}_{0, i + k - 2}(\eta^{(j - i + l - k + 4)/2})\Big) \\
& = \psi \tilde{p}_{k - 1, l - 1} \tilde{p}_{i - 1, j - 1} - \binom{i + k - 2}{i - 1} \tilde{p}_{i + k - 2, j + l} \\
& = \psi \tilde{p}_{i - 1, j - 1} \tilde{p}_{k - 1, l - 1} - \binom{i + k - 2}{i - 1} \tilde{p}_{i + k - 2, j + l}.
\end{split}
\end{equation}
On the other hand, since $f\big([J]\big) = 0$, we have
\begin{equation} \label{eq:deg1}
[f, \tilde{p}_{i, j}]\big([J]\big) = f(p_{i, j}) = p_{i - 2, j}.
\end{equation}
The identities \eqref{eq:deg2} and \eqref{eq:deg1} imply that for any polynomial $P$ in $\{p_{i, j}\}$ and $\psi$, we have
\begin{equation*}
f\Big(P\big(\{p_{i, j}\}, \psi\big)\Big) = \sD\Big(P\big(\{p_{i, j}\}, \psi\big)\Big),
\end{equation*}
where $\sD$ is the differential operator defined in \eqref{eq:defd} (see \cite{Pol07}, Section~3).

It remains to prove that $\psi \in \sT(J)$. To see this, we apply $\sD$ to the class $p_{1, 1}^2 \in \sT(J)$, which gives
\begin{equation*}
\sD(p_{1, 1}^2) = \psi p_{0, 0}^2 - \binom{0}{0} p_{0, 2} = g^2\psi - p_{0, 2}.
\end{equation*}
Hence $\psi = \big(\sD(p_{1, 1}^2) + p_{0, 2}\big)/g^2 \in \sT(J)$. \qed

\subsection{Corollary} --- \label{cor:tautisoms} \textit{For $i \geq 0$, there is an identity}
\begin{equation} \label{eq:basechange}
p_{0, 2i} = \pi^*\Bigg(\frac{1}{2^{i + 1}}\sum_{0 \leq j \leq i} \binom{i + 1}{j + 1} \psi^{i - j} \kappa_j + \psi^i\Bigg).
\end{equation}
\textit{Moreover, we have the following isomorphisms of $\bQ$-algebras \rleft similar to those in \textup{\eqref{eq:chowisoms}}\rright.}
\begin{equation} \label{eq:tautisoms}
\begin{tikzcd}[column sep=tiny,row sep=scriptsize]
\phantom{_g}\big(\bigoplus_{i = 0}^{d}\sT_{(0, 2i)}(J), \cdot\big) \arrow[leftarrow]{rr}{\sF}[swap]{\sim} & & \big(\bigoplus_{i = 0}^{d}\sT_{(2g, 2i)}(J), *\big) \arrow{dl}[pos=.61]{\ \pi_*}[swap]{\sim} \\
& \big(\sR(S), \cdot\big) \arrow{ul}{\pi^*}[swap]{\sim}
\end{tikzcd}
\end{equation}
\textit{In particular, the ring $\sR(S)$ may be regarded as a $\bQ$-subalgebra of $\big(\sT(J), \cdot\big)$ via $\pi^*$.}

\begin{proof}
By the isomorphisms \eqref{eq:chowisoms} we have
\begin{equation*}
p_{0, 2i} = \sF\big(\theta^{i + 1} \cdot [C]_{(2i)}\big) = \pi^*\pi_*\big(\theta^{i + 1} \cdot [C]_{(2i)}\big) = \pi^*\pi_*\big(\theta^{i + 1} \cdot [C]\big).
\end{equation*}
Hence it suffices to calculate $\pi_*\big(\theta^{i + 1} \cdot [C]\big)$. Applying \eqref{eq:pulltheta} and the projection formula, we get
\begin{align*}
\pi_*\big(\theta^{i + 1} \cdot [C]\big) & = p_*\Bigg(\bigg(\frac{1}{2}K + [x_0] + \frac{1}{2}\psi\bigg)^{i + 1}\Bigg) \\
& = \sum_{\substack{j + k + l = i + 1 \\ j, k, l \geq 0}} \frac{(i + 1)!}{j!k!l!}\frac{1}{2^{j + l}} p_*\big(K^j \cdot [x_0]^k \cdot \psi^l\big).
\end{align*}
Again by applying the projection formula to $p \colon C \to S$ and $x_0 \colon S \to C$, we find
\begin{equation*}
p_*\big(K^j \cdot [x_0]^k \cdot \psi^l\big) = \psi^l \cdot p_*\big(K^j \cdot [x_0]^k\big) = 
\begin{cases}
\psi^l \cdot \kappa_{j - 1} & \rif k = 0, \\
\psi^l \cdot x_0^*\big(K^j \cdot [x_0]^{k - 1}\big) = (-1)^{k - 1} \psi^i & \rif k \geq 1,
\end{cases}
\end{equation*}
with the convention $\kappa_{-1} = 0$. It follows that
\begin{align*}
\pi_*\big(\theta^{i + 1} \cdot [C]\big) & = \sum_{\substack{j + l = i + 1 \\ j, l \geq 0}} \frac{(i + 1)!}{j!l!}\frac{1}{2^{i + 1}} \psi^l \kappa_{j - 1} + \sum_{\substack{j + k + l = i + 1 \\ j, k, l \geq 0}} \frac{(i + 1)!}{j!k!l!}\frac{1}{2^{j + l}} (-1)^{k - 1} \psi^i \\
& \qquad - \sum_{\substack{j + l = i + 1 \\ j, l \geq 0}} \frac{(i + 1)!}{j!l!}\frac{1}{2^{i + 1}} (-1) \psi^i \\
& = \frac{1}{2^{i + 1}}\sum_{0 \leq j \leq i} \binom{i + 1}{j + 1} \psi^{i - j} \kappa_j + \bigg(\frac{1}{2} - 1 + \frac{1}{2}\bigg)^{i + 1} \psi^i + \bigg(\frac{1}{2} + \frac{1}{2}\bigg)^{i + 1} \psi^i \\
& = \frac{1}{2^{i + 1}}\sum_{0 \leq j \leq i} \binom{i + 1}{j + 1} \psi^{i - j} \kappa_j + \psi^i,
\end{align*}
which proves the identity \eqref{eq:basechange}.

Now since $\oplus_{i = 0}^{d}\sT_{(0, 2i)}(J)$ is generated by the classes $\{p_{0, 2i}\}$ and $\psi$, we have one inclusion $\oplus_{i = 0}^{d}\sT_{(0, 2i)}(J) \subset \pi^*\big(\sR(S)\big)$. For the other inclusion, it follows from \eqref{eq:basechange} and induction on~$j$ that one can also express $\pi^*(\kappa_j)$ as linear combinations of $\{p_{0, 2i}\}$ and $\psi$. So we have $\oplus_{i = 0}^{d}\sT_{(0, 2i)}(J) = \pi^*\big(\sR(S)\big)$, and the rest follows from \eqref{eq:chowisoms}.
\end{proof}

\secskip
\section{The Chow ring(s) of the relative infinite symmetric power}
\secskip

\noindent A nice way to treat the Chow rings of all $C^{[n]}$ simultaneously is to consider the relative infinite symmetric power $C^{[\infty]}$. It allows us to lift structures from the Chow ring of $J$. We recall the definition of $C^{[\infty]}$ and its Chow theories. Our reference is the paper of Moonen and Polishchuk \cite{MP10}, which generalizes the work of Kimura and Vistoli \cite{KV96} to the relative setting.

\subsectionempty We retain the notation in Sections~\ref{sec:notation} and~\ref{sec:notation2}. For $n \geq 1$, define maps $\varphi_n \colon C^{[n]} \to J$ and $\phi_n \coloneqq \varphi_n \circ \sigma_n \colon C^n \to J$, which send (locally) $n$ sections $x_1, \ldots, x_n$ of $C/S$ to the class $\sO_C(x_1 + \cdots + x_n - n x_0)$. Note that $\varphi_1 = \phi_1 = \iota \colon C \into J$. For convenience we set $\varphi_0 = \phi_0 \coloneqq o \colon S \to J$. We have the following diagram in addition to \eqref{eq:setting}.
\begin{equation} \label{eq:setting2}
\begin{tikzcd}
C^n \arrow{r}{\sigma_n} \arrow[bend left=30]{rr}{\phi_n} \arrow{rrd}[swap]{p^n} & C^{[n]} \arrow{r}{\varphi_n} \arrow{rd}[pos=.6]{p^{[n]}} & J \arrow{d}[pos=.4]{\pi} \\
& & S
\end{tikzcd}
\end{equation}

\subsectionempty We introduce the relative \textit{infinite symmetric power} of $C$, which is defined to be the ind-scheme
\begin{equation*}
C^{[\infty]} \coloneqq \varinjlim (S = C^{[0]} \into C \into C^{[2]} \into C^{[3]} \into \cdots).
\end{equation*}
Here the transition maps $\epsilon_n \colon C^{[n - 1]} \into C^{[n]}$ are given by adding a copy of $x_0$ (in particular $\epsilon_1 = x_0 \colon S \to C$). We write $p^{[\infty]} \colon C^{[\infty]} \to S$ for the limit of $p^{[n]} \colon C^{[n]} \to S$.

The collection of $\varphi_n \colon C^{[n]} \to J$ induces $\varphi \colon C^{[\infty]} \to J$. Moreover, the monoidal structure on $J$ can be lifted to $C^{[\infty]}$: the addition maps $\mu_{n, m} \colon C^{[n]} \times_S C^{[m]} \to C^{[n + m]}$ give rise to
\begin{equation*}
\mu \colon C^{[\infty]} \times_S C^{[\infty]} \to C^{[\infty]},
\end{equation*}
while the diagonal maps $\Delta_N \colon C^{[n]} \to C^{[Nn]}$ give rise to
\begin{equation*}
[N] \colon C^{[\infty]} \to C^{[\infty]} \rfor N \geq 0.
\end{equation*}

Unlike the case of $J$, there are (at least) two different notions of Chow ring for $C^{[\infty]}$. One is graded by codimension and the other by relative dimension. For this reason, we also distinguish the two gradings on $\CH(J)$: we write $\CH^\bullet(J) \coloneqq \big(\oplus_i \CH^i(J), \cdot\big)$ and $\CH_\bullet(J) \coloneqq \big(\oplus_i \CH_i(J), *\big)$.

\subsection{Definition} --- (i) \ The \textit{Chow cohomology} of $C^{[\infty]}$ of codimension $i$ is the inverse limit
\begin{equation*}
\CH^i(C^{[\infty]}) \coloneqq \varprojlim \big(\CH^i(S) \gets \CH^i(C) \gets \CH^i(C^{[2]}) \gets \CH^i(C^{[3]}) \gets \cdots\big),
\end{equation*}
where the transition maps are $\epsilon_n^* \colon \CH^i(C^{[n]}) \to \CH^i(C^{[n - 1]})$.

(ii) \ The \textit{Chow homology} of $C^{[\infty]}$ of relative dimension $i$ is the direct limit
\begin{equation*}
\CH_i(C^{[\infty]}) \coloneqq \varinjlim \big(\CH_i(S) \to \CH_i(C) \to \CH_i(C^{[2]}) \to \CH_i(C^{[3]}) \to \cdots\big),
\end{equation*}
where the transition maps are $\epsilon_{n, *} \colon \CH_i(C^{[n - 1]}) \to \CH_i(C^{[n]})$.

\secskip
By definition, an element in $\CH^i(C^{[\infty]})$ is a sequence $\alpha = (\alpha_n)_{n \geq 0}$ with $\alpha_n \in \CH^i(C^{[n]})$, such that $\epsilon_n^*(\alpha_n) = \alpha_{n - 1}$. An element in $\CH_i(C^{[\infty]})$ is represented by $\alpha \in \CH_i(C^{[n]})$ for some~$n$. We write $\CH^\bullet(C^{[\infty]}) \coloneqq \oplus_i \CH^i(C^{[\infty]})$ and $\CH_\bullet(C^{[\infty]}) \coloneqq \oplus_i \CH_i(C^{[\infty]})$. The Chow cohomology and homology of $C^{[\infty]} \times_S C^{[\infty]}$ are defined similarly.

\subsection{Remark} --- For $n \geq 1$, the transition map $\epsilon_n^*$ (\resp $\epsilon_{n, *}$) is surjective (\resp injective). To see this, one constructs a correspondence $\delta_n \colon C^{[n]} \vdash C^{[n - 1]}$ satisfying $\delta_n \circ \epsilon_n = \id$, so that $\epsilon_n^* \circ \delta_n^* = \id$ and $\delta_{n, *} \circ \epsilon_{n, *} = \id$. Let $\pr_i \colon C^n \to C$ be the $i$-th projection. Then consider the correspondence (see \cite{KV96}, Definition~1.8)
\begin{equation} \label{eq:defdelta}
\sum_{j = 0}^{n - 1} \sum_{1 \leq i_1 < \cdots < i_j \leq n}(-1)^{n - 1 - j}(\pr_{i_1}, \ldots, \pr_{i_j}, x_0, \ldots, x_0) \colon C^n \vdash C^{n - 1},
\end{equation}
which is $\fS_n$-invariant and descends to $\delta_n \colon C^{[n]} \vdash C^{[n - 1]}$. One verifies that $\delta_n$ has the desired property. In particular, the map $\CH^\bullet(C^{[\infty]}) \to \CH(C^{[n]})$ (\resp $\CH(C^{[n]}) \to \CH_\bullet(C^{[\infty]})$) is surjective (\resp injective).

\subsectionempty Both $\CH^\bullet(C^{[\infty]})$ and $\CH_\bullet(C^{[\infty]})$ are equipped with a ring structure: there is the intersection product ``$\cdot$'' on $\CH^\bullet(C^{[\infty]})$ and the \textit{Pontryagin product} ``$*$'' on $\CH_\bullet(C^{[\infty]})$, the latter defined in the same way as \eqref{eq:pontryagin}.

Note that unlike the case of $\CH(J)$, one cannot define both ring structures on the same object $\CH^\bullet(C^{[\infty]})$ or $\CH_\bullet(C^{[\infty]})$. We do, however, have a \textit{cap product}
\begin{equation*}
\CH^\bullet(C^{[\infty]}) \times \CH_\bullet(C^{[\infty]}) \xto{\cap} \CH_\bullet(C^{[\infty]}),
\end{equation*}
which sends $\alpha = (\alpha_n)_{n \geq 0} \in \CH^i(C^{[\infty]})$ and $\beta \in \CH_j(C^{[m]}) \subset \CH_j(C^{[\infty]})$ to $\alpha \cap \beta \coloneqq \alpha_m \cdot \beta \in \CH_{j - i}(C^{[m]}) \subset \CH_{j - i}(C^{[\infty]})$. One verifies that ``$\cap$'' is well-defined.

For $N \geq 0$, the multiplication $[N] \colon C^{[\infty]} \to C^{[\infty]}$ induces
\begin{equation*}
[N]^* \colon \CH^\bullet(C^{[\infty]}) \to \CH^\bullet(C^{[\infty]}), \rand [N]_* \colon \CH_\bullet(C^{[\infty]}) \to \CH_\bullet(C^{[\infty]}).
\end{equation*}
Further, the map $\varphi \colon C^{[\infty]} \to J$ induces morphisms of $\bQ$-algebras
\begin{equation*}
\varphi^* \colon \CH^\bullet(J) \to \CH^\bullet(C^{[\infty]}), \rand \varphi_* \colon \CH_\bullet(C^{[\infty]}) \to \CH_\bullet(J),
\end{equation*}
which allow us to relate the Chow theories of $C^{[\infty]}$ and $J$. We begin with $\CH^\bullet(C^{[\infty]})$.

\subsection{Chow cohomology} --- Recall the class $\psi \in \CH^1(S)$. Again we keep the same notation $\psi$ for its pull-backs to schemes over $S$. For $n \geq 1$, define
\begin{equation} \label{eq:defxi}
\sO_{C^{[n]}}(1) \coloneqq \sO_{C^{[n]}}\big(\epsilon_n(C^{[n - 1]}) + n \psi\big),
\end{equation}
and denote by $\xi_n \in \CH^1(C^{[n]})$ the first Chern class of $\sO_{C^{[n]}}(1)$. We set $\xi_0 = 0$. Then we have $\epsilon_n^*(\xi_n) = \xi_{n - 1}$, which yields a class
\begin{equation*}
\xi \coloneqq (\xi_n)_{n \geq 0} \in \CH^1(C^{[\infty]}).
\end{equation*}

There is an alternative description of $\sO_{C^{[n]}}(1)$ (essentially due to Schwarzenberger \cite{Sch63}). Let $\sL$ be the pull-back of the Poincaré line bundle $\sP$ via $\iota \times_S \id_J \colon C \times_S J \to J \times_S J$. For $n \geq 0$, define the sheaf
\begin{equation*}
E_n \coloneqq \pr_{2, *}\Big(\pr_1^*\big(\sO_C(n x_0)\big) \otimes \sL\Big),
\end{equation*}
where $\pr_1 \colon C \times_S J \to C$ and $\pr_2 \colon C \times_S J \to J$ are the two projections. There is a canonical isomorphism $C^{[n]} \simeq \bP(E_n)$, under which $\sO_{C^{[n]}}(1)$ corresponds to the line bundle $\sO_{\bP(E_n)}(1)$.

When $n \geq 2g - 1$, the sheaf $E_n$ is locally free over $J$ (\ie a vector bundle). Then we have isomorphisms of $\bQ$-algebras
\begin{equation} \label{eq:chowbundle}
\CH(C^{[n]}) \simeq \CH\big(\bP(E_n)\big) \simeq \CH(J)[\xi_n]\left/\big<P(\xi_n)\big>,\right.
\end{equation}
where $P(\xi_n) \coloneqq \sum_{i = 0}^{n - g + 1} c_{n - g + 1 - i}(E_n) \xi_n^i$ is a polynomial in $\xi_n$ of degree $n - g + 1$ with coefficients in $\CH(J)$.

\subsectionempty The following results are merely reinterpretations of \eqref{eq:chowbundle} (\cite{MP10}, Theorem~1.4).

(i) \ When $n \geq \max\{2g, i + g + 1\}$, the transition map $\epsilon_n^* \colon \CH^i(C^{[n]}) \to \CH^i(C^{[n - 1]})$ becomes an isomorphism.

(ii) \ The map $\varphi^* \colon \CH^\bullet(J) \to \CH^\bullet(C^{[\infty]})$ induces an isomorphism of $\bQ$-algebras
\begin{equation} \label{eq:isophi}
\Phi \colon \CH^\bullet(J)[t] \xto{\sim} \CH^\bullet(C^{[\infty]}),
\end{equation}
which sends $\alpha \in \CH^\bullet(J)$ to $\varphi^*(\alpha)$ and $t$ to the class $\xi$.

\subsection{Chow homology} --- We switch to $\CH_\bullet(C^{[\infty]})$. First, choose an integer $n \geq 2g + d$ (recall that $d = \dim(S)$), and we identify $C^{[n]}$ with $\bP(E_n)$. Consider the short exact sequence
\begin{equation*}
0 \to \sO_{C^{[n]}}(-1) \to \varphi_n^*(E_n) \to Q \to 0,
\end{equation*}
where $Q$ is the universal quotient bundle of $\varphi_n^*(E_n)$, and is of rank $n - g$. Define
\begin{equation*}
\Gamma \coloneqq c_{n - g}(Q) \in \CH_g(C^{[n]}) \subset \CH_g(C^{[\infty]}).
\end{equation*}

Next, choose $n \geq 2g + d + 1$. Define
\begin{equation*}
L \coloneqq \varphi_n^*\big([o]\big) \cdot c_{n - g - 1}(Q) \in \CH_1(C^{[n]}) \subset \CH_1(C^{[\infty]}),
\end{equation*}
where $[o] \coloneqq \big[o(S)\big] \in \CH_0(J)$ is the class of the zero section. One can show that both classes $\Gamma \in \CH_g(C^{[\infty]})$ and $L \in \CH_1(C^{[\infty]})$ are independent of $n$.

\subsectionempty Here are the corresponding results for $\CH_\bullet(C^{[\infty]})$ (\cite{MP10}, Theorem~1.11).

(i) \ The map $s \colon \CH_\bullet(J) \to \CH_\bullet(C^{[\infty]})$ given by
\begin{equation*}
s(\alpha) \coloneqq \varphi^*(\alpha) \cap \Gamma
\end{equation*}
is a section of $\varphi_* \colon \CH_\bullet(C^{[\infty]}) \to \CH_\bullet(J)$. It respects the Pontryagin products on both sides.

(ii) \ The section $s$ induces an isomorphism of $\bQ$-algebras
\begin{equation} \label{eq:isopsi}
\Psi \colon \CH_\bullet(J)[t] \xto{\sim} \CH_\bullet(C^{[\infty]}),
\end{equation}
which sends $\alpha \in \CH_\bullet(J)$ to $s(\alpha)$ and $t$ to the class $L$.

(iii) \ Under the isomorphism $\Psi$, the push-forward $\varphi_*$ is the evaluation at zero, and the action of $\xi \cap$ is the derivation $\rd/\rd t$.

\subsection{Remark} --- There are explicit expressions for $\Gamma$ and $L$ (see \cite{MP10}, Corollary~1.13)
\begin{equation} \label{eq:formulagamma}
\Gamma = \frac{1}{g!(N - 1)^g} \vast(\frac{\log\Big(1 + \psi \cdot [N]_*\big([C]\big)\Big) - N \log\big(1 + \psi \cdot [C]\big)}{N\psi}\vast)^g,
\end{equation}
with $N \geq 2$, and
\begin{equation} \label{eq:formulal}
L = \frac{\log\big(1 + \psi \cdot [C]\big) - \log\bigg(1 + \psi \cdot s\Big(\big[\iota(C)\big]\Big)\bigg)}{\psi}.
\end{equation}
Here we distinguish $[C] \in \CH_1(C) \subset \CH_1(C^{[\infty]})$ from $\big[\iota(C)\big] \in \CH_1(J)$.

\subsection{Fourier transform} --- By \cite{KV96}, Theorems~3.13 and~3.18, the Fourier transform on $J$ can be lifted to $C^{[\infty]}$. Recall the class 
$\ell \in \CH^1(J \times_S J)$. Define
\begin{align*}
\ell_{\infty, \infty} & \coloneqq (\varphi \times_S \varphi)^*(\ell) \in \CH^1(C^{[\infty]} \times_S C^{[\infty]}), \\
\xi \times_S \xi & \coloneqq \pr_1^*(\xi) \cdot \pr_2^*(\xi) \in \CH^2(C^{[\infty]} \times_S C^{[\infty]}),
\end{align*}
where $\pr_1, \pr_2 \colon C^{[\infty]} \times_S C^{[\infty]} \to C^{[\infty]}$ are the two projections. Then the expression
\begin{equation} \label{eq:corf}
\exp(\ell_{\infty, \infty} + \xi \times_S \xi)
\end{equation}
is an \textit{upper correspondence} in the sense of \loccit, Definition~3.2. It induces an isomorphism of $\bQ$-algebras
\begin{equation*}
\sF \colon \CH_\bullet(C^{[\infty]}) \xto{\sim} \CH^\bullet(C^{[\infty]}).
\end{equation*}
We have $\sF(L) = \xi$, and the following commutative diagram.
\begin{equation*}
\begin{tikzcd}
\CH_\bullet(J) \arrow{r}{\sF}[swap]{\sim} \arrow{d}[swap]{s} & \CH^\bullet(J) \arrow{d}{\varphi^*} \\
\CH_\bullet(C^{[\infty]}) \arrow{r}{\sim}[swap]{\sF} & \CH^\bullet(C^{[\infty]})
\end{tikzcd}
\end{equation*}
Further, the inverse $\sF^{-1} \colon \CH^\bullet(C^{[\infty]}) \xto{\sim} \CH_\bullet(C^{[\infty]})$ is given by the \textit{lower correspondence} in the sense of \loccit, Definition~3.17
\begin{equation} \label{eq:corfinv}
(-1)^g \exp(-\ell_{\infty, \infty}) \cap \Big(\exp^*\big((L \times_S L) * (\Gamma \times_S \Gamma)\big)\Big).
\end{equation}
where $\exp^*$ means the exponential power series with respect to ``$*$''.

\subsection{Remark} --- The class $\ell_{\infty, \infty}$ can be described somewhat explicitly. For $n, m \geq 0$, define
\begin{equation*}
\ell_{n, m} \coloneqq (\varphi_n \times_S \varphi_m)^*(\ell) \in \CH^1(C^{[n]} \times_S C^{[m]}).
\end{equation*}
We have $\ell_{n, 0} = \ell_{0, m} = 0$, and we have seen in \eqref{eq:pulll} that
\begin{equation*}
\ell_{1, 1} = [\Delta] - \pr_1^*\big([x_0]\big) - \pr_2^*\big([x_0]\big) - \psi \rin \CH^1(C \times_S C).
\end{equation*}
For $1 \leq i \leq n$, let $\pr_i \colon C^n \to C$ be the $i$-th projection. Then for $n, m \geq 1$, there is an identity (see \cite{KV96}, Proposition-Definition~3.10)
\begin{equation} \label{eq:pulll2}
(\phi_n \times_S \phi_m)^*(\ell) = (\sigma_n \times_S \sigma_m)^*(\ell_{n, m}) = \sum_{i = 1}^n \sum_{j = 1}^m (\pr_i \times_S \pr_j)^*(\ell_{1, 1}),
\end{equation}
which holds in $\CH^1(C^n \times_S C^m)^{\fS_n \times \fS_m}$.

\secskip
\section{The tautological ring(s) of the relative infinite symmetric power}
\secskip

\noindent We define the tautological cohomology $\sR^\bullet(C^{[\infty]})$ and homology $\sR_\bullet(C^{[\infty]})$ of the relative infinite symmetric power. We prove a tautological analogue of the isomorphisms $\Phi$ and $\Psi$ in the previous section, that both $\sR^\bullet(C^{[\infty]})$ and $\sR_\bullet(C^{[\infty]})$ are polynomial algebras over $\sT(J)$. The connections between the two notions of tautological ring are thus established. Throughout this section we work in the setting of \eqref{eq:setting} and \eqref{eq:setting2}.

\subsectionempty First observe that the rings $\sR(C^{[n]})$ are stable under pull-backs and push-forwards via the maps $\epsilon_n \colon C^{[n - 1]} \into C^{[n]}$. In fact, the maps $\epsilon_n$ can be lifted to $\id_{C^{n - 1}} \times_S x_0 \colon C^{n - 1} \to C^n$. Then for $\alpha \in \sR(C^n)$ and $\beta \in \sR(C^{n - 1})$, we have
\begin{align*}
(\id_{C^{n - 1}} \times_S x_0)^*(\alpha) & = \pr_{1, \ldots, n - 1, *}\Big(\alpha \cdot \pr_n^*\big([x_0]\big)\Big) \in \sR(C^{n - 1}), \\
(\id_{C^{n - 1}} \times_S x_0)_*(\beta) & = \pr_{1, \ldots, n - 1}^*(\beta) \cdot \pr_n^*\big([x_0]\big) \in \sR(C^n).
\end{align*}
Here $\pr_{1, \ldots, n - 1}$ (\resp $\pr_n$) is the projection of $C^n$ to the first $n - 1$ factors (\resp $n$-th factor), and is tautological in the sense of Remark~\ref{rem:deffp}. 

The stability of $\sR(C^{[n]})$ under $\epsilon_n^*$ and $\epsilon_{n, *}$ allows us to pass to $C^{[\infty]}$.

\subsection{Definition} --- (i) \ The \textit{tautological cohomology} of $C^{[\infty]}$ is the inverse limit
\begin{equation*}
\sR^\bullet(C^{[\infty]}) \coloneqq \varprojlim \big(\sR(S) \gets \sR(C) \gets \sR(C^{[2]}) \gets \sR(C^{[3]}) \gets \cdots\big),
\end{equation*}
where the transition maps are $\epsilon_n^* \colon \sR(C^{[n]}) \to \sR(C^{[n - 1]})$.

(ii) \ The \textit{tautological homology} of $C^{[\infty]}$ is the direct limit
\begin{equation*}
\sR_\bullet(C^{[\infty]}) \coloneqq \varinjlim \big(\sR(S) \to \sR(C) \to \sR(C^{[2]}) \to \sR(C^{[3]}) \to \cdots\big),
\end{equation*}
where the transition maps are $\epsilon_{n, *} \colon \sR(C^{[n - 1]}) \to \sR(C^{[n]})$.

\secskip
Note that $\sR^\bullet(C^{[\infty]})$ (\resp $\sR_\bullet(C^{[\infty]})$) inherits a grading from $\CH^\bullet(C^{[\infty]})$ (\resp $\CH_\bullet(C^{[\infty]})$). It is immediate that $\sR^\bullet(C^{[\infty]})$ is stable under the intersection product ``$\cdot$''. Since the addition map $\mu_{n, m} \colon C^{[n]} \times_S C^{[m]} \to C^{[n + m]}$ lifts to the identity (tautological) map $C^n \times_S C^m \to C^{n + m}$, we also know that $\sR_\bullet(C^{[\infty]})$ is stable under the Pontryagin product ``$*$''. It follows that $\sR^\bullet(C^{[\infty]})$ (\resp $\sR_\bullet(C^{[\infty]})$) is a graded $\bQ$-subalgebra of $\CH^\bullet(C^{[\infty]})$ (\resp $\CH_\bullet(C^{[\infty]})$).

We list a few properties of $\sR^\bullet(C^{[\infty]})$ and $\sR_\bullet(C^{[\infty]})$, which partly reveal their links with the ring $\sT(J)$.

\subsection{Proposition} \label{prop:proptaut} --- (i) \ \textit{We have $\xi \in \sR^1(C^{[\infty]})$, $\Gamma \in \sR_g(C^{[\infty]})$ and $L \in \sR_1(C^{[\infty]})$.}

(ii) \ \textit{The ring $\sR^\bullet(C^{[\infty]})$ \rleft \resp $\sR_\bullet(C^{[\infty]})$\rright\ is stable under $[N]^*$ \rleft \resp $[N]_*$\rright, for all $N \geq 0$.}

(iii) \ \textit{The Fourier transform $\sF$ induces an isomorphism}
\begin{equation*}
\sF \colon \sR_\bullet(C^{[\infty]}) \xto{\sim} \sR^\bullet(C^{[\infty]}).
\end{equation*}

(iv) \ \textit{The cap product restricts to a map}
\begin{equation*}
\sR^\bullet(C^{[\infty]}) \times \sR_\bullet(C^{[\infty]}) \xto{\cap} \sR_\bullet(C^{[\infty]}).
\end{equation*}

\begin{proof}
Statement~(ii) follows from the fact that the diagonal map $\Delta_N \colon C^{[n]} \to C^{[Nn]}$ lifts to $C^n \to C^{Nn}$, which is tautological. Statement~(iv) is straightforward.

For (i), by \eqref{eq:defxi} and \eqref{eq:formulagamma} we have $\xi \in \sR^1(C^{[\infty]})$ and $\Gamma \in \sR_g(C^{[\infty]})$. Moreover by \eqref{eq:formulal}, to prove that $L \in \sR_1(C^{[\infty]})$ it suffices to show that $s\Big(\big[\iota(C)\big]\Big) \in \sR_1(C^{[\infty]})$. This is further reduced to proving that $\varphi^*\Big(\big[\iota(C)\big]\Big) \in \sR^{g - 1}(C^{[\infty]})$ by the definition of the section $s$.

In fact, we can prove for any $\alpha \in \sT(J)$ that $\varphi^*(\alpha) \in \sR^\bullet(C^{[\infty]})$. First by Theorem~\ref{thm:structure}, we know that $\big(\sT(J), \cdot\big)$ is generated by the classes $\{p_{i, j}\}$ and $\psi$. Since $\varphi^*(\psi) \in \sR^1(C^{[\infty]})$, it remains to prove that $\varphi^*(p_{i, j}) \in \sR^\bullet(C^{[\infty]})$ for all possible $i$ and $j$. Here we can actually calculate the pull-back of $p_{i, j}$ via $\phi_n = \varphi_n \circ \sigma_n \colon C^n \to J$, for all $n \geq 0$. The procedure is similar to that of Lemma~\ref{lem:pulltheta}: we chase through the following cartesian squares.
\begin{equation*}
\begin{tikzcd}
C \times_S C^n \arrow{r}{\iota \times_S \id_{C^n}} \arrow{d}[swap]{\id_C \times_S \phi_n} & J \times_S C^n \arrow{r}{\pr_2} \arrow{d}{\id_J \times_S \phi_n} & C^n \arrow{d}{\phi_n} \\
C \times_S J \arrow{r}[swap]{\iota \times_S \id_J} \arrow{d}[swap]{\pr_1} & J \times_S J \arrow{r}[swap]{\pr_2} \arrow{d}{\pr_1} & J \\
C \arrow{r}[swap]{\iota} & J & \phantom{C \times_S C}
\end{tikzcd}
\end{equation*}
Then we find
\begin{equation} \label{eq:pullpij}
\phantom{(5.1)}\phi_n^*\bigg(\sF\Big(\theta^{(j - i + 2)/2} \cdot \big[\iota(C)\big]\Big)\bigg) = \pr_{2, *}\Big(\pr_1^*\big(\iota^*(\theta)^{(j - i + 2)/2}\big) \cdot \exp\big((\iota \times_S \phi_n)^*(\ell)\big)\Big),
\end{equation}
where $\pr_1 \colon C \times_S C^n \to C$ and $\pr_2 \colon C \times_S C^n \to C^n$ are the two projections. By definition $\phi_n^*(p_{i, j})$ is just the codimension $(i + j)/2$ component of the right-hand side of \eqref{eq:pullpij}. Further by \eqref{eq:pulltheta} and \eqref{eq:pulll2}, we have explicit expressions for $\iota^*(\theta)$ and $(\iota \times_S \phi_n)^*(\ell)$ in terms of tautological classes. It follows that $\phi_n^*(p_{i, j}) \in \sR(C^n)$, and hence $\varphi^*(p_{i, j}) \in \sR^\bullet(C^{[\infty]})$.

Finally to prove (iii), we observe that the correspondences in \eqref{eq:corf} and \eqref{eq:corfinv} that define $\sF$ and $\sF^{-1}$ only involve tautological classes. 
\end{proof}

Now we state and prove the main result of this section. To be coherent, we write $\sT^\bullet(J) \coloneqq \big(\oplus_i \sT^i(J), \cdot\big)$ and $\sT_\bullet(J) \coloneqq \big(\oplus_i \sT_i(J), *\big)$.

\subsection{Theorem} \label{thm:poltaut} --- \textit{The isomorphisms $\Phi$ and $\Psi$ in \textup{\eqref{eq:isophi}} and \textup{\eqref{eq:isopsi}} restrict to isomorphisms of $\bQ$-algebras}
\begin{align}
\Phi|_{\sT^\bullet(J)[t]} \colon \sT^\bullet(J)[t] & \xto{\sim} \sR^\bullet(C^{[\infty]}), \label{eq:tautisophi} \\
\Psi|_{\sT_\bullet(J)[t]} \colon \sT_\bullet(J)[t] & \xto{\sim} \sR_\bullet(C^{[\infty]}). \label{eq:tautisopsi}
\end{align}

\secskip
The plan is to prove \eqref{eq:tautisopsi} first, and then deduce \eqref{eq:tautisophi} by Fourier duality. We begin with an elementary lemma.

\subsection{Lemma} \label{lem:trivial} --- \textit{Let $A$ be a commutative $\bQ$-algebra, and $B$ be a $\bQ$-subalgebra of the polynomial algebra $A[t]$. Assume that $t \in B$, and that $B$ is stable under derivation $\rd/\rd t$. Then we have}
\begin{equation*}
B = \ev(B)[t],
\end{equation*}
\textit{where $\ev \colon A[t] \to A$ is the evaluation at zero.}

\begin{proof}
Take an element $P(t) = b_0 + b_1t + \cdots + b_nt^n$ in $B$. Since $B$ is stable under derivation, we have $(\rd/\rd t)^n\big(P(t)\big) = n!b_n \in B$, so that $b_n \in B$. Then since $t \in B$, we have $b_nt^n \in B$ and $P(t) - b_nt^n \in B$. By induction, we find that all coefficients $b_i$ are in $B$. It follows that $\ev(B)[t] \subset B$. On the other hand, we know that $b_i = \ev\Big((\rd/\rd t)^i\big(P(t)/i!\big)\Big)$ with $(\rd/\rd t)^i\big(P(t)/i!\big) \in B$. Hence $b_i \in \ev(B)$, which proves the other inclusion $B \subset \ev(B)[t]$.
\end{proof}

Consider the push-forward map $\varphi_* \colon \CH_\bullet(C^{[\infty]}) \to \CH_\bullet(J)$ which, under $\Psi$, corresponds to the evaluation at zero. The proof of the following proposition is a bit involved and relies essentially on Theorem~\ref{thm:structure} and Corollary~\ref{cor:tautisoms}.

\subsection{Proposition} \label{prop:pushtaut} --- \textit{We have $\varphi_*\big(\sR_\bullet(C^{[\infty]})\big) = \sT_\bullet(J)$.}

\begin{proof}
By Theorem~\ref{thm:structure}, we know that $\sT_\bullet(J)$ is generated by $\Big\{\theta^{(j - i + 2)/2} \cdot \big[\iota(C)\big]_{(j)}\Big\}$ and $o_*(\psi)$ (recall that $o = \varphi_0$ is the zero section). Consider the class $\eta = K/2 + [x_0] + \psi/2 \in \sR(C)$, which by \eqref{eq:pulltheta} is equal to $\iota^*(\theta)$. Then we have
\begin{equation*}
\iota_*(\eta^{(j - i + 2)/2}) = \theta^{(j - i + 2)/2} \cdot \big[\iota(C)\big],
\end{equation*}
so that $\theta^{(j - i + 2)/2} \cdot \big[\iota(C)\big]$ is in the image $\varphi_*\big(\sR_\bullet(C^{[\infty]})\big)$ (recall that $\iota = \varphi_1$). Moreover, we have shown in Proposition~\ref{prop:proptaut}~(ii) that $\sR_\bullet(C^{[\infty]})$ is stable under $[N]_*$ for all $N \geq 0$. This implies that the components $\theta^{(j - i + 2)/2} \cdot \big[\iota(C)\big]_{(j)}$ are also in the image $\varphi_*\big(\sR_\bullet(C^{[\infty]})\big)$. Since all generators of $\sT_\bullet(J)$ are in $\varphi_*\big(\sR_\bullet(C^{[\infty]})\big)$, we obtain the inclusion $\sT_\bullet(J) \subset \varphi_*\big(\sR_\bullet(C^{[\infty]})\big)$.

To prove the other inclusion, observe that $\varphi_*\big(\sR_\bullet(C^{[\infty]})\big)$, being the union of $\varphi_{n, *}\big(\sR(C^{[n]})\big)$ for $n \geq 0$, is also the union of $\phi_{n, *}\big(\sR(C^n)\big)$ for $n \geq 0$. Then it is enough to prove that  $\phi_{n, *}\big(\sR(C^n)\big) \subset \sT(J)$ for all $n \geq 0$. This is done by an explicit calculation in terms of the generators of $\sR(C^n)$.

As $\phi_{n, *}$ factors through $\phi_{n + 1, *}$, we may assume $n \geq 2$. The ring $\sR(C^n)$ is then generated by $\{\kappa_i\}$, $\psi$, $\{K_j\}$ and $\big\{[x_{0, j}]\big\}$, and $\big\{[\Delta_{k, l}]\big\}$. We make a change of variables
\begin{equation*}
\eta_j \coloneqq \frac{1}{2}K_j + [x_{0, j}] + \frac{1}{2}\psi,
\end{equation*}
so that $\sR(C^n)$ is also generated by $\{\kappa_i\}$, $\psi$, $\{\eta_j\}$ and $\big\{[x_{0, j}]\big\}$, and $\big\{[\Delta_{k, l}]\big\}$. Let $\alpha \in \sR(C^n)$ be a monomial in those generators. We would like to show that $\phi_{n, *}(\alpha) \in \sT(J)$. 

A first step is to separate the variables $\{\kappa_i\}$ and $\psi$ from the rest. Write $\alpha = \beta \cdot \gamma$, with $\beta$ collecting all factors of $\{\kappa_i\}$ and $\psi$. Then $\beta$ is the pull-back of a class $\beta_0 \in \sR(S)$ via $p^n \colon C^n \to S$. Since $p^n = \pi \circ \phi_n$ (recall that $\pi \colon J \to S$), we find
\begin{equation*}
\phi_{n, *}(\alpha) = \phi_{n, *}\big(p^{n, *}(\beta_0) \cdot \gamma\big) = \phi_{n, *}\big(\phi_n^*\pi^*(\beta_0) \cdot \gamma\big) = \pi^*(\beta_0) \cdot \phi_{n, *}(\gamma).
\end{equation*}
Thanks to the isomorphisms \eqref{eq:tautisoms}, we have $\pi^*(\beta_0) \in \sT(J)$. So it remains to prove that $\phi_{n, *}(\gamma) \in \sT(J)$, or in other words, we may assume that $\alpha$ is a monomial in $\{\eta_j\}$, $\big\{[x_{0, j}]\big\}$ and $\big\{[\Delta_{k, l}]\big\}$ only.

A second step is to eliminate multiplicities in the variables $\big\{[\Delta_{k, l}]\big\}$. Consider for example $[\Delta] = [\Delta_{1, 2}] \in \sR(C^2)$. Denote by $\Delta \colon C \to C^2$ the diagonal map, and by $\pr_2 \colon C^2 \to C$ the second projection. Then we have
\begin{equation*}
[\Delta]^2 = \Delta_*\Big(\Delta^*\big([\Delta]\big)\Big) = -\Delta_*(K) = -\Delta_*\big(\Delta^*\pr_2^*(K)\big) = -[\Delta] \cdot K_2.
\end{equation*}
By pulling back to $C^n$, we obtain for $1 \leq k < l \leq n$
\begin{equation*}
[\Delta_{k, l}]^2 = -[\Delta_{k, l}] \cdot K_l = -[\Delta_{k, l}] \cdot \big(2\eta_l - 2[x_{0, l}] - \psi\big).
\end{equation*}
Together with the first step, this allows us to reduce to the case where $\alpha$ is a monomial in $\{\eta_j\}$, $\big\{[x_{0, j}]\big\}$ and $\big\{[\Delta_{k, l}]\big\}$, with multiplicity at most $1$ for each $[\Delta_{k, l}]$.

Further, we may permute the indices of the $\big\{[\Delta_{k, l}]\big\}$ factors by applying the identity
\begin{equation} \label{eq:permute}
[\Delta_{k, l}] \cdot [\Delta_{l, m}] = [\Delta_{k, m}] \cdot [\Delta_{l, m}].
\end{equation}
More precisely, if $I = \{i_1, i_2, \ldots, i_q\}$ is a subset of $\{1, \ldots, n\}$, we define the symbol
\begin{equation*}
[\Delta_I] \coloneqq [\Delta_{i_1, i_2}] \cdot [\Delta_{i_1, i_3}] \cdots [\Delta_{i_1, i_q}].
\end{equation*}
It follows from \eqref{eq:permute} that $[\Delta_I]$ is well-defined. Also we have identities $[\Delta_I] \cdot \eta_{i_1} = [\Delta_I] \cdot \eta_{i_2} = \cdots = [\Delta_I] \cdot \eta_{i_q}$ and $[\Delta_I] \cdot [x_{0, i_1}] = [\Delta_I] \cdot [x_{0, i_2}] = \cdots = [\Delta_I] \cdot [x_{0, i_q}]$. So for $r, s \geq 0$, we can write
\begin{equation*}
(\eta^r[x_0]^s)_{\Delta_I} \coloneqq [\Delta_I] \cdot \eta_{i_1}^r \cdot [x_{0, i_1}]^s.
\end{equation*}

Combining with the first two steps, we may assume that the $\big\{[\Delta_{k, l}]\big\}$ factors of $\alpha$ take the~form
\begin{equation*}
[\Delta_{I_1}] \cdot [\Delta_{I_2}] \cdots [\Delta_{I_m}],
\end{equation*}
where the $I_k$ are subsets of $\{1, 2, \ldots, n\}$ satisfying $I_k \cap I_l = \emptyset$ for $k \neq l$. This means we are reduced to the case where $\alpha$ is of the form
\begin{equation*}
\alpha = \prod_{k \in \{1, \ldots, m\}} (\eta^{r_k}[x_0]^{s_k})_{\Delta_{I_k}} \cdot \prod_{\substack{j \in \{1, \ldots, n\} \\ j \notin \cup_k I_k}} (\eta_j^{u_j}[x_{0, j}]^{v_j}),
\end{equation*}
with $I_1, \ldots, I_m$ and $\{j\}$ pairwise disjoint. In this case, the calculation of $\phi_{n, *}(\alpha)$ is rather straightforward: it follows almost from the definitions of $[N]$ and ``$*$'' that
\begin{equation*}
\phi_{n, *}(\alpha) = \sideset{}{^*}\prod_{k \in \{1, \ldots, m\}} [\#I_k]_*\iota_*(\eta^{r_k}[x_0]^{s_k}) * \sideset{}{^*}\prod_{\substack{j \in \{1, \ldots, n\} \\ j \notin \cup_k I_k}} \iota_*(\eta^{u_j}[x_0]^{v_j}),
\end{equation*}
where $\prod^*$ stands for product with respect to ``$*$'', and $\#I_k$ the cardinality of $I_k$.

Now since $\sT(J)$ is stable under $[N]_*$ and ``$*$'', the last step is to prove that $\iota_*\big(\eta^r[x_0]^s\big) \in \sT(J)$ for all $r, s \geq 0$. By the identity $\eta = \iota^*(\theta)$, we have $\iota_*\big(\eta^r[x_0]^s\big) = \theta^r \cdot \iota_*\big([x_0]^s\big)$, which further reduces to showing that $\iota_*\big([x_0]^s\big) \in \sT(J)$. Note that $\iota_*\big([x_0]^0\big) = \big[\iota(C)\big] \in \sT(J)$ and $\iota_*\big([x_0]\big) = [o] \in \sT(J)$ (recall that $[o]$ is the class of the zero section). For $s \geq 2$, we have
\begin{equation*}
\iota_*\big([x_0]^s\big) = \iota_*x_{0, *}x_0^*\big([x_0]^{s - 1}\big) = o_*\big((-\psi)^{s - 1}\big) \in \sT(J).
\end{equation*}
The proof of the inclusion $\varphi_*\big(\sR_\bullet(C^{[\infty]})\big) \subset \sT_\bullet(J)$ is thus completed.
\end{proof}

\subsection{Proof of Theorem~\textup{\ref{thm:poltaut}}} --- By Proposition~\ref{prop:proptaut}~(i) and~(iv), we know that $L \in \sR_\bullet(C^{[\infty]})$, and that $\sR_\bullet(C^{[\infty]})$ is stable under $\xi \cap -$. Then the isomorphism \eqref{eq:tautisopsi} follows immediately from Lemma~\ref{lem:trivial} and Proposition~\ref{prop:pushtaut}. By applying the Fourier transform $\sF$ and by Proposition~\ref{prop:proptaut}~(iii), we also obtain \eqref{eq:tautisophi}. \qed

\subsection{Remark} --- Previously, Moonen and Polishchuk considered much bigger tautological rings of $C^{[\infty]}$ and $J$, for which they obtained similar results as Theorem~\ref{thm:poltaut} (\cite{MP10}, Corollary~8.6). The advantage of our minimalist version, is that one can use the $\fsl_2$-machinery on~$\sT(J)$ to study enumerative problems on $\sR(C^{[n]})$. 

\secskip
\section{Tautological relations and Gorenstein properties}
\secskip

\noindent The $\fsl_2$-action on the Jacobian side provides relations between tautological classes. Using these relations, we study the Gorenstein property for $\sR(\sM_{g, 1})$ (\resp $\sR(\sM_g)$). We also formulate the corresponding Gorenstein property for the tautological ring $\sT(\sJ_{g, 1})$ of the universal Jacobian. Then we prove that $\sT(\sJ_{g, 1})$ being Gorenstein is equivalent to $\sR(\sC_{g, 1}^{[n]})$ being Gorenstein for all $n \geq 0$. Computation confirms the Gorenstein properties for small $g$, and suggests when these properties may not hold.

\subsection{Relations via \texorpdfstring{$\fsl_2$}{sl2}} --- We explain how the $\fsl_2$-action gives relations in $\sT(J)$. Following the isomorphisms \eqref{eq:tautisoms}, we identify $\sR(S)$ with $\oplus_{i = 0}^d \sT_{(0, 2i)}(J)$ via the map $\pi^*$. Then we obtain relations in $\sR(S)$ by restriction.

The idea is due to Polishchuk \cite{Pol05}. By Theorem~\ref{thm:structure}, the space of polynomial relations between $\{p_{i, j}\}$ and $\psi$ is stable under the action of $\sD$. In other words, if $P$ is a polynomial in $\{p_{i, j}\}$ and $\psi$, then $P\big(\{p_{i, j}\}, \psi\big) = 0$ implies $\sD\Big(P\big(\{p_{i, j}\}, \psi\big)\Big) = 0$. Now consider monomials
\begin{equation*}
\alpha = \psi^s p_{i_1, j_1}^{r_1} p_{i_2, j_2}^{r_2} \cdots p_{i_m, j_m}^{r_m} \rwith I \coloneqq r_1 i_1 + r_2 i_2 + \cdots + r_m i_m > 2g.
\end{equation*}
By definition $\alpha \in \oplus_j\CH_{(I, j)}(J)$. But since $I > 2g$, we know from the decomposition \eqref{eq:beauville2} that $\alpha = 0$. In terms of the Dutch house, the class $\alpha$ is simply outside the house. It follows that we have relations
\begin{equation*}
\alpha = 0, \ \ \sD(\alpha) = 0, \ \ \sD^2(\alpha) = 0, \ \ \ldots
\end{equation*}

\subsectionempty This argument leads to the following formal definition. Let $i, j$ run through all integers such that $i \leq j + 2$ and that $i + j$ is even. Define
\begin{equation*}
\sA \coloneqq \bQ\big[\{x_{i, j}\}, y\big]\left/\big<x_{0, 0} - g, \{x_{i, j}\}_{i < 0}, \{x_{i, j}\}_{j < 0}, \{x_{i, j}\}_{j > 2g - 2}\big>.\right.
\end{equation*}
In other words, the ring $\sA$ is a polynomial ring in variables $\{x_{i, j}\}$ and $y$, with the convention that $x_{0, 0} = g$ and $x_{i, j} = 0$ for $i < 0 $ or $j < 0$ or $j > 2g - 2$ (same as the classes $\{p_{i, j}\}$). We introduce a bigrading $\sA = \oplus_{i, j}\sA_{(i, j)}$ by the requirement that $x_{i, j} \in \sA_{(i, j)}$ and $y \in \sA_{(0, 2)}$. Define operators $E, F$ and $H$ on $\sA$ by
\begin{align*}
E \colon \sA_{(i, j)} & \to \sA_{(i + 2, j)} & \alpha & \mapsto x_{2, 0} \cdot \alpha, \\
F \colon \sA_{(i, j)} & \to \sA_{(i - 2, j)} & \alpha & \mapsto F(\alpha), \\
H \colon \sA_{(i, j)} & \to \sA_{(i, j)} & \alpha & \mapsto (i - g)\alpha,
\end{align*}
where 
\begin{equation*}
\begin{split}
F \coloneqq \frac{1}{2}\sum_{i, j, k, l}\Bigg(y x_{i - 1, j - 1} x_{k - 1, l - 1} & - \binom{i + k - 2}{i - 1} x_{i + k - 2, j + l}\Bigg) \partial x_{i, j} \partial x_{k, l} \\
& + \sum_{i, j}x_{i - 2, j}\partial x_{i, j}.
\end{split}
\end{equation*}
It is not difficult to verify that the operators above generate a $\bQ$-linear representation $\fsl_2 \to \End_\bQ(\sA)$. Theorem~\ref{thm:structure} can then be reformulated as the existence of a surjective morphism of $\fsl_2$-representations $\sA \to \sT(J)$, which maps $x_{i, j}$ to $p_{i, j}$ and $y$ to $\psi$.

\subsectionempty Denote by $\Mon_{(i, j)}$ the set of non-zero monomials in $\sA_{(i, j)}$ (also without $x_{0, 0}$ as a factor). For convenience we set $\Mon_{(0, 0)} \coloneqq \{1\}$. Then consider the quotient
\begin{equation} \label{eq:ttilde}
\tilde{\sT} \coloneqq \sA\left/\Big<\big\{F^\nu(\Mon_{(i, j)})\big\}_{i > 2g, \nu \geq 0}\Big>.\right.
\end{equation}
The ring $\tilde{\sT}$ inherits a bigrading $\tilde{\sT} = \oplus_{i, j}\tilde{\sT}_{(i, j)}$ from $\sA$. The operators $E, F$ and $H$ induce operators on $\tilde{\sT}$, which we denote by $e, f$ and $h$. Again we obtain a representation $\fsl_2 \to \End_\bQ(\tilde{\sT})$. Further, since $e^{g + 1} = f^{g + 1} = 0$, we formally define the Fourier transform on $\tilde{\sT}$ by
\begin{equation*}
\sF \coloneqq \exp(e) \exp(-f) \exp(e).
\end{equation*}

\begin{figure}
\centering
\includegraphics[height=.5\textheight]{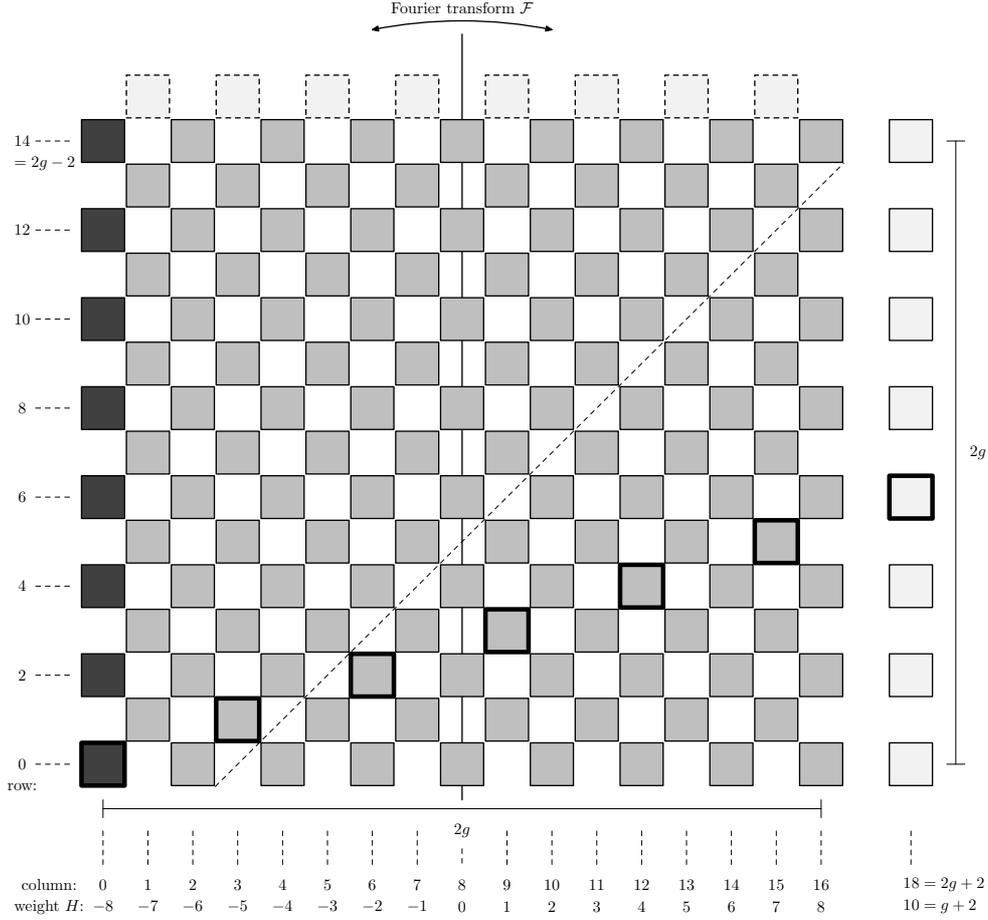}
\caption{Producing relations ($g = 8$).}
\label{fig:relations}
\end{figure}

We also define the subring $\tilde{\sR} = \oplus_i \tilde{\sT}_{(0, 2i)}$, graded by $\tilde{\sR} = \oplus_i \tilde{\sR}^i$ with $\tilde{\sR}^i \coloneqq \tilde{\sT}_{(0, 2i)}$. We~have
\begin{align*}
\tilde{\sR}^i & = \sA_{(0, 2i)}\left/\Big<\big\{F^I(\Mon_{(2I, 2i)})\big\}_{I > g}\Big>\right.\\
& = \sA_{(0, 2i)}\left/\big<F^{g + 1}(\Mon_{(2g + 2, 2i)})\big>.\right.
\end{align*}
Figure~\ref{fig:relations} illustrates the construction of $\tilde{\sR}$: take monomials on the $(2g + 2)$-th column of the Dutch house (white blocks), and then apply $g + 1$ times the operator $F$ to obtain relations between the generators (black blocks).

To summarize this formal approach, we have the following proposition.

\subsection{Proposition} --- \textit{There are surjective maps}
\begin{equation} \label{eq:bigtheta}
\Theta \colon \tilde{\sT} \onto \sT(J), \rand \Theta|_{\tilde{\sR}} \colon \tilde{\sR} \onto \sR(S),
\end{equation}
\textit{which map $x_{i, j}$ to $p_{i, j}$ and $y$ to $\psi$}. \qed

\secskip
We look at $\tilde{\sR}$ in more details. The following lemma shows one can eliminate certain monomials that produce trivial relations.

\subsection{Lemma} --- \textit{For all $\alpha \in \Mon_{(2g + 2, 2i)}$ of the form $\alpha = x_{2, 0} \cdot \beta$, we have $F^{g + 1}(\alpha) = 0$.}

\begin{proof}
We know that $F^{g + 1}(\alpha) = F^{g + 1}(x_{2, 0} \cdot \beta) = F^{g + 1}E(\beta)$. Using a combinatorial identity for $\fsl_2$-representations (see for example \cite{Moo09}, Lemma~2.4), we find
\begin{equation*}
F^{g + 1}E(\beta) = EF^{g + 1}(\beta) - (g + 1)(g + 1 - g - 1 + 1 - 1)F^g(\beta) = EF^{g + 1}(\beta).
\end{equation*}
On the other hand $F^{g + 1}(\beta) \in \sA_{(-2, 2i)} = 0$, which implies $F^{g + 1}(\alpha) = 0$.
\end{proof}

As a result, if we write $\mon_{(2g + 2, 2i)} \subset \Mon_{(2g + 2, 2i)}$ for the subset of monomials without $x_{2, 0}$ as a factor, then we have
\begin{equation} \label{eq:rtilde}
\tilde{\sR}^i = \sA_{(0, 2i)}\left/\big<F^{g + 1}(\mon_{(2g + 2, 2i)})\big>.\right. 
\end{equation}
The bold blocks in Figure~\ref{fig:relations} describe the lower bound of $i$ such that $\mon_{(2g + 2, 2i)}$ is non-empty. In numerical terms, we have $x_{3, 1}^{2i} \in \sA_{(6i, 2i)}$ and that $6i > 2g$ means $i > g/3$. Hence $\mon_{(2g + 2, 2i)} = \emptyset$ for all $i \leq \lfloor g/3 \rfloor$.

\subsection{Proposition} --- \textit{The elements $x_{0, 2}, \ldots, x_{0, 2 \lfloor g/3 \rfloor}$ and $y$ generate $\tilde{\sR}$, with no relations in~$\tilde{\sR}^i$ for $i \leq \lfloor g/3 \rfloor$.}

\secskip
Together with \eqref{eq:basechange}, this gives a new proof of Ionel's result (\cite{Ion05}, Theorem~1.5), that $\sR(S)$ is generated by $\kappa_1, \ldots, \kappa_{\lfloor g/3 \rfloor}$ and $\psi$.

\begin{proof}
The second part is immediate after \eqref{eq:rtilde} and the fact that $\mon_{(2g + 2, 2i)} = \emptyset$ for all $i \leq \lfloor g/3 \rfloor$. For the first part, the goal is to relate all $x_{0, 2i}$ with $i > g/3$ to the elements $x_{0, 2}, \ldots, x_{0, 2 \lfloor g/3 \rfloor}$ and $y$, and the idea is to use specific monomials to get these relations.

We proceed by induction. Suppose all $\{x_{0, 2j}\}_{g/3 < j < i}$ can be expressed in terms of the elements $x_{0, 2}, \ldots, x_{0, 2 \lfloor g/3 \rfloor}$ and $y$. Then consider the monomial $x_{3, 1}^{2i} \in \sA_{(6i, 2i)}$. Applying $3i$-times the operator~$F$ we get $F^{3i}(x_{3, 1}^{2i}) \in \sA_{(0, 2i)}$, which vanishes in $\tilde{\sR}$. On the other hand, by going through the definition of $F$ we find
\begin{equation*}
F^{3i}(x_{3, 1}^{2i}) = c x_{0, 2i} + \alpha,
\end{equation*}
where $\alpha$ is a polynomial in $\{x_{0, 2j}\}_{j < i}$ and $y$. It remains to show that $c$ is non-zero.

The observation is the following: when we apply the operator $F$, the minus sign occurs every time two factors ($x_{i, j}, x_{k, l}$) are merged into one ($x_{i + k - 2, j + l}$). If we start from $x_{3, 1}^{2i}$ and arrive at $x_{0, 2i}$, no matter how we proceed we have to do the merging $(2i - 1)$-times. This means all non-zero summands of $c$ are of the form $(-1)^{2g - 1}$ times a positive number, hence negative. Then the sum $c$ is also negative.
\end{proof}

\subsection{Computing \texorpdfstring{$\sR(\sM_{g, 1})$}{R(M_g,1)}} --- Our colleague Li Ma made a C$++$ program that computes $\tilde{\sR}$ for a given genus $g$. It calculates relations and outputs the dimension of each component $\tilde{\sR}^i$.

Meanwhile, based on an algorithm developed by Liu and Xu \cite{LX12}, Bergvall computed the intersection numbers in $\sR(\sM_{g, 1})$ for many values of $g$ (see \cite{Ber11}, Section~4.2). It then gives the dimensions of the Gorenstein quotient $\sG(\sM_{g, 1})$, which is the quotient of $\sR(\sM_{g, 1})$ by those classes that pair zero with all opposite degree classes. Note that this computation is formal and does not involve actual relations in $\sR(\sM_{g, 1})$.

There are surjective maps $\tilde{\sR} \onto \sR(\sM_{g, 1}) \onto \sG(\sM_{g, 1})$. Our computation shows that for $g \leq 19$, the dimensions of $\tilde{\sR}$ and $\sG(\sM_{g, 1})$ are equal, which means we have $\tilde{\sR} \simeq \sR(\sM_{g, 1}) \simeq \sG(\sM_{g, 1})$. In particular, we can confirm the following (for $g \leq 9$ this has been obtained independently by Bergvall; see \cite{Ber11}, Section~4.4).

\subsection{Corollary} --- \textit{The ring $\sR(\sM_{g, 1})$ is Gorenstein for $g \leq 19$.} \qed

\secskip
However, the computer output is negative for $g = 20$ and some greater values of $g$. There the dimensions of $\tilde{\sR}$ are simply not symmetric. Again by comparing with the dimensions of $\sG(\sM_{g, 1})$, we know exactly how many relations are missing. The numbers are listed in Table~\ref{tab:missing} below. Note that for $g \geq 25$, we only calculated a range near the middle degree, presuming that $\sR(\sM_{g, 1})$ behaves well near the top.

\begin{table}
\tabcolsep = 10pt
\centering
{\renewcommand*\arraystretch{1.2}
\begin{tabular}{c|r@{: }l|r@{: }l}
\hline
$g$ & \multicolumn{2}{c|}{$\sM_{g, 1}$} & \multicolumn{2}{c}{$\phantom{\Big(}\sM_g\phantom{\Big)}$} \\
\hline
$\leq 19$ & \multicolumn{2}{c|}{OK} & \multicolumn{2}{c}{OK} \\
$20$ & codim $10$ & $1$ missing & \multicolumn{2}{c}{OK} \\
$21$ & codim $11$ & $1$ missing & \multicolumn{2}{c}{OK} \\
$22$ & codim $11$ & $1$ missing & \multicolumn{2}{c}{OK} \\
$23$ & codim $12$ & $3$ missing & \multicolumn{2}{c}{OK} \\
$24$ & codim $13$ & $2$ missing & codim $12$ & $1$ missing \\
& codim $12$ & $4$ missing \\
$25$ & codim $13$ & $5$ missing & codim $12$ & $1$ missing \\
& codim $12$ & $1$ missing \\
$26$ & codim $14$ & $6$ missing & codim $13$ & $1$ missing \\
& codim $13$ & $6$ missing \\
$27$ & codim $15$ & $3$ missing & codim $14$ & $1$ missing \\
& codim $14$ & $11$ missing & codim $13$ & $1$ missing \\
& codim $13$ & $1$ missing \\
$28$ & codim $15$ & $10$ missing & codim $14$ & $2$ missing \\
& codim $14$ & $10$ missing \\
\hline
\end{tabular}}
\vspace{\baselineskip}
\caption{Computer output for $g \leq 28$.}
\label{tab:missing}
\end{table}

\subsectionempty We include a brief discussion about $\sR(\sM_g)$. Recall that $\sR(\sM_g)$ is the $\bQ$-subalgebra of $\CH(\sM_g)$ generated by~$\{\kappa_i\}$. Faber's original conjecture predicts that $\sR(\sM_g)$ is Gorenstein with socle in degree~$g - 2$ (\cite{Fab99}, Conjecture~1).

Denote by $q \colon \sM_{g, 1} \to \sM_g$ the forgetful map. It is not difficult to verify that $q_*\big(\sR(\sM_{g, 1})\big) = \sR(\sM_g)$. More precisely, for $\psi^s \kappa_{i_1}^{r_1} \cdots \kappa_{i_m}^{r_m} \in \sR(\sM_{g, 1})$ we have
\begin{equation*}
q_*(\psi^s \kappa_{i_1}^{r_1} \cdots \kappa_{i_m}^{r_m}) = \kappa_{s - 1} \kappa_{i_1}^{r_1} \cdots \kappa_{i_m}^{r_m},
\end{equation*}
with the convention that $\kappa_{-1} = 0$.

Together with \eqref{eq:basechange}, the identity above allows us to push relations in $\sR(\sM_{g, 1})$ forward to $\sR(\sM_g)$. We used another computer program to do the work. Then we recover the well-known result of Faber and Zagier, that $\sR(\sM_g)$ is Gorenstein for $g \leq 23$. Note that when $20 \leq g \leq 23$, the missing relations in $\sR(\sM_{g, 1})$ do not affect the Gorenstein property for $\sR(\sM_g)$.

From $g = 24$ on, the computer output is again negative; see Table~\ref{tab:missing}. Our computation for $g \leq 28$ suggests that we obtain exactly the same set of relations as the Faber-Zagier relations (see \cite{PP13}, Section~0.2). Notably in the crucial case of $g = 24$, we have not found the missing relation in codimension~$12$. It is not known whether in theory we obtain the same relations.

\subsectionempty We continue to study the Gorenstein properties for $\sR(\sC_{g, 1}^{[n]})$. A key step is to establish a similar property for the Jacobian side. Denote by $\sJ_{g, 1}$ the universal Jacobian over $\sM_{g, 1}$. The tautological ring $\sT(\sJ_{g, 1})$ is thus defined.

The following lemma locates the expected socle for $\sT(\sJ_{g, 1})$.

\subsection{Lemma} --- \textit{We have $\sT^i(\sJ_{g, 1}) = 0$ for $i > 2g - 1$, and}
\begin{equation*}
\sT^{2g - 1}(\sJ_{g, 1}) = \sT^{2g - 1}_{(2g - 2)}(\sJ_{g, 1}) \simeq \bQ.
\end{equation*}

\begin{proof}
The surjective map $\phi_g \colon \sC_{g, 1}^g \onto \sJ_{g, 1}$ induces $\phi_g^* \colon \CH(\sJ_{g, 1}) \into \CH(\sC_{g, 1}^g)$ which, by~\eqref{eq:tautisophi}, restricts to an injective map $\phi_g^* \colon \sT(\sJ_{g, 1}) \into \sR(\sC_{g, 1}^g)$. Then it follows from \eqref{eq:socle} that $\sT^i(\sJ_{g, 1}) = 0$ for $i > 2g - 1$, and that $\sT^{2g - 1}(\sJ_{g, 1})$ is at most $1$-dimensional.. 

Recall from \eqref{eq:tautisoms} that $\sT^{g - 1}_{(2g - 2)}(\sJ_{g, 1}) \simeq \sR^{g - 1}(\sM_{g, 1}) \simeq \bQ$. Applying the Fourier transform, we obtain
\begin{equation*}
\sT^{2g - 1}_{(2g - 2)}(\sJ_{g, 1}) = \sF\big(\sT^{g - 1}_{(2g - 2)}(\sJ_{g, 1})\big) \simeq \bQ.
\end{equation*}
So $\sT^{2g - 1}(\sJ_{g, 1})$ is indeed $1$-dimensional, and is concentrated in $\sT^{2g - 1}_{(2g - 2)}(\sJ_{g, 1})$.
\end{proof}

With the socle condition verified, we consider for $0 \leq i \leq 2g - 1$ the pairing
\begin{equation} \label{eq:pairjac}
\sT^i(\sJ_{g, 1}) \times \sT^{2g - 1 - i}(\sJ_{g, 1}) \xto{\cdot} \sT^{2g - 1}(\sJ_{g, 1}) \simeq \bQ.
\end{equation}
There is the following analogue of Speculation~\ref{spec:faberpower}.

\subsection{Speculation} --- \textit{For $0 \leq i \leq 2g - 1$, the pairing \eqref{eq:pairjac} is perfect. In other words, the ring $\sT(\sJ_{g, 1})$ is Gorenstein with socle in degree $2g - 1$.}

\secskip
The Dutch house gives an nice interpretation. In Figure~\ref{fig:pairjac}, the socle component is located precisely in the upper right corner. Assuming the Gorenstein property, one would expect a rotational symmetry about the center of the picture. Together with the reflection symmetry about the middle vertical line (given by the Fourier transform), it would then imply a mysterious reflection symmetry about the middle horizontal line. In particular, one should have $\sT_{(j)}(\sJ_{g, 1}) = 0$ for $j > 2g - 2$, and using the grading in \eqref{eq:beauville2}, a one-to-one correspondence between $\sT_{(i, j)}(\sJ_{g, 1})$ and $\sT_{(i, 2g - 2 - j)}(\sJ_{g, 1})$.

\begin{figure}
\centering
\includegraphics[height=.5\textheight]{Pics_TaleTaut_4.mps}
\caption{Pairings in $\sT(\sJ_{g, 1})$ ($g = 8$).}
\label{fig:pairjac}
\end{figure}

\subsection{Computing \texorpdfstring{$\sT(\sJ_{g, 1})$}{T(J_g,1)}} --- Once again with the help of Li Ma, we computed the ring~$\tilde{\sT}$ defined in \eqref{eq:ttilde} and its pairings. The computer output shows that for $g \leq 7$, we do have $\tilde{\sT}^{2g - 1} \simeq \bQ$ and perfect pairings between $\tilde{\sT}^i$ and $\tilde{\sT}^{2g - 1 - i}$. Then since the surjective map $\tilde{\sT} \onto \sT(\sJ_{g, 1})$ is an isomorphism at the socle level, it is in fact an isomorphism. In particular, we can confirm the following.

\subsection{Corollary} \label{cor:calcttilde} --- \textit{The ring $\sT(\sJ_{g, 1})$ is Gorenstein for $g \leq 7$.} \qed

\secskip
For $g = 8$ and some greater values of $g$, however, the relations we find are not sufficient to deduce the Gorenstein property for $\sT(\sJ_{g, 1})$.

Further, one of the main results of this section is an equivalence of Gorenstein properties.

\subsection{Theorem} \label{thm:equiv} --- \textit{Fix $g > 0$. The following three statements are equivalent:}

(i) \ \textit{The ring $\sT(\sJ_{g, 1})$ is Gorenstein;}

(ii) \ \textit{The ring $\sR(\sC_{g, 1}^{[n]})$ is Gorenstein for some $n \geq 2g - 1$;}

(iii) \ \textit{The ring $\sR(\sC_{g, 1}^{[n]})$ is Gorenstein for all $n \geq 0$.}

\begin{proof}
First assume $n \geq 2g - 1$. Recall that $\sC_{g, 1}^{[n]}$ is a $\bP^{n - g}$-bundle over $\sJ_{g, 1}$. Then we obtain from \eqref{eq:chowbundle} and \eqref{eq:tautisophi} an isomorphism of $\bQ$-algebras
\begin{equation*}
\sR(\sC_{g, 1}^{[n]}) \simeq \sT(\sJ_{g, 1})[\xi_n]\left/\big<P(\xi_n)\big>,\right.
\end{equation*}
where $P(\xi_n)$ is a polynomial in $\xi_n$ of degree $n - g + 1$ with coefficients in $\sT(\sJ_{g, 1})$ (one shows by induction that the coefficients of $P(\xi_n)$ are in $\varphi_{n, *}\big(\sR(\sC_{g, 1}^{[n]})\big)$, which by Proposition~\ref{prop:pushtaut} equals $\sT(\sJ_{g, 1})$). In particular, we have for the socle components
\begin{equation*}
\sR^{g - 1 + n}(\sC_{g, 1}^{[n]}) \simeq \sT^{2g - 1}(\sJ_{g, 1}) \cdot \xi_n^{n - g} \simeq \bQ.
\end{equation*}

For $0 \leq i \leq g - 1 + n$, we write $\sR^i(\sC_{g, 1}^{[n]}) \simeq \oplus_j \sT^{i - j}(\sJ_{g, 1}) \cdot \xi_n^j$ with $\max\{0, i - 2g + 1\} \leq j \leq \min\{i, n - g\}$. Then the pairing \eqref{eq:pairsym} corresponds to
\begin{equation*}
\Big(\bigoplus_j \sT^{i - j}(\sJ_{g, 1}) \cdot \xi_n^j\Big) \times \Big(\bigoplus_k \sT^{g - 1 + n - i - k}(\sJ_{g, 1}) \cdot \xi_n^k\Big) \xto{\cdot} \sT^{2g - 1}(\sJ_{g, 1}) \cdot \xi_n^{n - g} \simeq \bQ.
\end{equation*}
On the other hand, observe that
\begin{equation*}
\big(\sT^{i - j}(\sJ_{g, 1}) \cdot \xi_n^j\big) \cdot \big(\sT^{g - 1 + n - i - k}(\sJ_{g, 1}) \cdot \xi_n^k\big) = 0 \rif j + k < n - g.
\end{equation*}
In other words, if we choose a suitable basis for each $\sR^i(\sC_{g, 1}^{[n]})$, then the pairing matrix of~\eqref{eq:pairsym} is block triangular. Moreover, the blocks on the diagonal correspond exactly to the case $j + k = n - g$, \ie the pairing
\begin{equation*}
\big(\sT^{i - j}(\sJ_{g, 1}) \cdot \xi_n^j\big) \times \big(\sT^{2g - 1 - i + j}(\sJ_{g, 1}) \cdot \xi_n^{n - g - j}\big) \xto{\cdot} \sT^{2g - 1}(\sJ_{g, 1}) \cdot \xi_n^{n - g} \simeq \bQ.
\end{equation*}
which in turn corresponds to the pairing
\begin{equation*}
\sT^{i - j}(\sJ_{g, 1}) \times \sT^{2g - 1 - i + j}(\sJ_{g, 1}) \xto{\cdot} \sT^{2g - 1}(\sJ_{g, 1}) \simeq \bQ.
\end{equation*}
In total, saying that \eqref{eq:pairsym} is perfect for all $0 \leq i \leq g - 1 + n$, is equivalent to saying that \eqref{eq:pairjac} is perfect for all $0 \leq i \leq 2g - 1$. This settles the proof for (some or all) $n \geq 2g - 1$.

For the remaining $n$, the observation is that $\sR(\sC_{g, 1}^{[n]})$ being Gorenstein implies $\sR(\sC_{g, 1}^{[n - 1]})$ being Gorenstein. To see this, take a class $\alpha \in \sR^i(\sC_{g, 1}^{[n - 1]})$ that pairs zero with all $\beta \in \sR^{g + n - 2 - i}(\sC_{g, 1}^{[n - 1]})$. Recall the correspondence $\delta_n \colon \sC_{g, 1}^{[n]} \vdash \sC_{g, 1}^{[n - 1]}$ defined in \eqref{eq:defdelta}, and consider the pairing of $\delta_n^*(\alpha)$ with all $\gamma \in \sR^{g + n - 1 - i}(\sC_{g, 1}^{[n]})$. It then follows from the projection formula that each summand of $\delta_n^*(\alpha) \cdot \gamma$ is zero, and hence $\delta_n^*(\alpha) \cdot \gamma = 0$. Assuming $\sR(\sC_{g, 1}^{[n]})$ Gorenstein, we get $\delta_n^*(\alpha) = 0$. Since $\delta_n^*$ is injective, we have $\alpha = 0$ and $\sR(\sC_{g, 1}^{[n - 1]})$ Gorenstein. The proof is thus completed.
\end{proof}

By Corollary~\ref{cor:calcttilde} and the discussion in Section~\ref{sec:unpoint}, we obtain an immediate consequence.

\subsection{Corollary} --- \textit{The ring $\sR(\sC_{g, 1}^{[n]})$ \rleft\resp $\sR(\sC_g^{[n]})$\rright\ is Gorenstein for $g \leq 7$ and for all $n \geq 0$.} \qed

\secskip
Finally, we conclude by proposing an alternative description of the tautological rings. As is explained in the introduction, this new description is coherent with Polishchuk's philosophy in \cite{Pol05}, and appears more geometric than the Gorenstein expectation.

\subsection{Conjecture} \label{conj:motive} --- \textit{The map $\Theta \colon \tilde{\sT} \to \sT(\sJ_{g, 1})$ in \textup{\eqref{eq:bigtheta}} is an isomorphism.}

\secskip

\end{document}